\documentclass{amsart}

\newtheorem{theorem}{Theorem}[section]
\newtheorem{lemma}[theorem]{Lemma}
\newtheorem{proposition}[theorem]{Proposition}
\newtheorem{question}[theorem]{Question}
\theoremstyle{definition}
\newtheorem{definition}[theorem]{Definition}
\theoremstyle{remark}
\newtheorem{remark}[theorem]{Remark}
\newtheorem{corollary}[theorem]{Corollary}

\usepackage{empheq}
\usepackage{tikz}
\usepackage{amssymb}
\usepackage[shortlabels]{enumitem}
\usepackage{bussproofs}
\usepackage{multicol}

\usepackage[utf8]{inputenc}

\usepackage[left=2.0cm,%
                right=2.0cm,%
                top=2.5cm,%
                bottom=2.5cm,%
                headheight=12pt,%
                a4paper]{geometry}%
                
\newcommand{\Wedge}[2]{\ensuremath{\sideset{}{^{#1}}\bigwedge_{#2}}}
\newcommand{\Vee}[2]{\ensuremath{\sideset{}{^{#1}}\bigvee_{#2}}}

\begin{document}

\title{On Infinitary G\"odel logics}

\author{Nicholas Pischke}
\address{Hoch-Weiseler Str. 46, Butzbach, 35510, Hesse, Germany}
\email{pischkenicholas@gmail.com}
\date{\today}


\begin{abstract}
We study propositional and first-order G\"odel logics over infinitary languages which are motivated semantically by corresponding interpretations into the unit interval $[0,1]$. We provide infinitary Hilbert-style calculi for the particular (propositional and first-order) cases with con-/disjunctions of countable length and prove corresponding completeness theorems by extending the usual Lindenbaum-Tarski construction to the infinitary case for a respective algebraic semantics via complete linear Heyting algebras. We provide infinitary hypersequent calculi and prove corresponding cut-elimination theorems in the Sch\"utte-Tait-style. Initial observations are made regarding truth-value sets other than $[0,1]$.
\end{abstract}

\maketitle

\section{Introduction}
Infinitary logics in a classical setting go back to \cite{Kar1964,ST1958,Tar1958} and over time became influential in various areas of mathematical logic like (finite) model theory, set theory and also formal arithmetic, among others. Model theoretically, they pose an interesting challenge since the usual propositional and first-order properties (like e.g. compactness) become more intertwined with set-theoretic principles (like e.g. large cardinal axioms).

We study infinitary extensions of G\"odel logics. In the finitary setting, G\"odel logics arose historically from a sequence of propositional finite-valued logics given by G\"odel \cite{Goe1932} to show that intuitionistic logic does not have a finite characteristic matrix. These were extended to an infinite-valued variant by Dummett \cite{Dum1959} and the whole collection is today especially studied in the context of intermediate logics. Further, G\"odel logics have been characterized as one of three main instances of t-norm based fuzzy logics by H\'ajek \cite{Haj1998}. First-order versions were first described by Horn \cite{Hor1969} and later rediscovered by Takeuti and Titani \cite{TT1984} under the name of intuitionistic fuzzy logics. The infinitary versions studied here assume a similar position among both the infinitary intermediate and infinitary fuzzy logics and the present work is thus, in that way, also a particular case study of these classes. For that purpose, G\"odel logics pose an especially interesting case as they, in the intermediate context, are logics with many classical properties but which, at the same time, are distinct enough from classical or intuitionistic logic to still pose interesting methodical challenges for the adaption of well-known results (e.g. like interpolation). An example of this phenomenon in the infinitary setting is the work \cite{Agu2016} by Aguilera where he studies analogues of the compactness results for classical infinitary logics in a G\"odel setting where, although the classical results stay true modulo appropriate reformulations, the methods of Skolem functions used classically had to be reformulated using a certain theory of fuzzy ultraproducts.

By now, the work \cite{Agu2016} is the only paper on infinitary G\"odel logics and many interesting problems arising from generalizations of the classical case have remained open, like the study of propositional variants, the development of infinitary calculi and appropriate completeness theorems in the propositional and first-order case as well as infinitary structural proof theory, among others. We study all these previously named topics (in variable depth) and in particular prove the relevant completeness and cut-elimination theorems (where one naturally restricts to the instances with conjunctions and disjunctions of countable length and finitary quantifiers like in the classical (see \cite{Kar1964}) and intuitionistic cases (see \cite{Nad1978})). At the end, we consider truth-value sets different than $[0,1]$ and extend the results to some of these cases as well.
\section{Propositional Infinitary G\"odel Logics}
\subsection{Syntax and Fragments}
Let $\kappa$ be any cardinal number. The \emph{infinitary propositional language associated with $\kappa$} is given by 
\[
\mathcal{L}_\kappa:\phi ::= \bot\mid x\mid (\phi\rightarrow\phi)\mid (\phi\land\phi)\mid(\phi\lor\phi)\mid\bigwedge\Phi\mid\bigvee\Phi
\]
where we have $x\in Var_\kappa:=\{x_\lambda\mid\lambda\in\kappa\}$ and $\Phi$ is a set of formulas of size $<\kappa$. For the other classical operators, we define
\begin{enumerate}
\item $\neg\phi:=\phi\rightarrow\bot$,
\item $\top:=\neg\bot$,
\item $\phi\leftrightarrow\psi:=(\phi\rightarrow\psi)\land(\psi\rightarrow\phi)$.
\end{enumerate}
Given some formula $\phi$, we denote the set of all subformulas (including $\phi$) by $\mathrm{sub}(\phi)$ and the set of variables of $\phi$ by $\mathrm{var}(\phi)$. Both of these naturally extend to sets $\Gamma$.

We mostly deal with the special case of $\kappa=\omega_1$ and in that context, we will mostly write $\bigwedge_{i\in\omega}\phi_i$ for $\bigwedge\{\phi_i\mid i\in\omega\}$ and $\bigvee_{i\in\omega}\phi_i$ for $\bigvee\{\phi_i\mid i\in\omega\}$ where $(\phi_i)_{i\in\omega}$ is a countable family of formulas.

Related to that particular instance of $\kappa=\omega_1$, we will also need the notion of a \emph{fragment}. These fragments are (possibly countable) sublanguages of $\mathcal{L}_{\omega_1}$ for which a Lindenbaum-Tarski construction is, nevertheless, still possible and they form a cornerstone of the proof of the completeness theorem. While these fragments, in particular the notation $\mathcal{L}_A$, originate from the connection of (classical) infinitary logic with admissible sets in the sense of Barwise \cite{Bar1969}, we only need and use the following syntactic definition, in similarity to Nadel \cite{Nad1978} in the context of infinitary intuitionistic logic.
\begin{definition}
A \emph{(distributive) fragment of $\mathcal{L}_{\omega_1}$} is a set $\mathcal{L}_A\subseteq\mathcal{L}_{\omega_1}$ such that
\begin{enumerate}
\item $\bot\in\mathcal{L}_A$,
\item $\phi\in\mathcal{L}_A$ implies $\mathrm{sub}(\phi)\subseteq\mathcal{L}_A$,
\item $\phi,\psi\in\mathcal{L}_A$ implies $\phi\circ\psi\in\mathcal{L}_A$ for $\circ\in\{\land,\lor,\rightarrow\}$,
\item $\phi,\bigwedge_{i\in\omega}\phi_i\in\mathcal{L}_A$ implies $\bigwedge_{i\in\omega}(\phi\rightarrow\phi_i)\in\mathcal{L}_A$,
\item $\phi,\bigvee_{i\in\omega}\phi_i\in\mathcal{L}_A$ implies $\bigwedge_{i\in\omega}(\phi_i\rightarrow\phi)\in\mathcal{L}_A$,
\item $\phi,\bigwedge_{i\in\omega}\phi_i\in\mathcal{L}_A$ implies $\bigwedge_{i\in\omega}(\phi\lor\phi_i)\in\mathcal{L}_A$.
\end{enumerate}
\end{definition}
The important kind of fragments will be countable ones. In particular, we will consider the smallest fragments containing some set of formulas.
\begin{lemma}\label{lem:genfrag}
For any $\Gamma\subseteq\mathcal{L}_{\omega_1}$, there is a smallest (w.r.t. $\subseteq$) distributive fragment $\mathrm{frag}(\Gamma)$ such that $\Gamma\subseteq\mathrm{frag}(\Gamma)$. If $\Gamma$ is countable, then $\mathrm{frag}(\Gamma)$ is also countable. 
\end{lemma}
\subsection{The Standard Semantics and $\mathsf{G}_\kappa$}
We now introduce the standard semantics for the language $\mathcal{L}_\kappa$ and the resulting logics of semantic consequence. This standard semantics naturally extends the usual finitary standard semantics for propositional G\"odel logics. For G\"odel logics, being one of the prime examples of many-valued logics, the first important parameter in that context is that of the truth-value set. We fix this to be $[0,1]$ for the major part of the paper and discuss other choices only in part at the end.
\begin{definition}\label{def:standardpropsem}
An \emph{$\mathcal{L}_\kappa$-G\"odel-evaluation} is a function $v:\mathcal{L}_\kappa\to [0,1]$ such that
\begin{enumerate}
\item $v(\bot)=0$,
\item $v(\phi\land\psi)=\min\{v(\phi),v(\psi)\}$,
\item $v(\phi\lor\psi)=\max\{v(\phi),v(\psi)\}$,
\item $v(\phi\rightarrow\psi)=v(\phi)\Rightarrow v(\psi)$ where $x\Rightarrow y:=\begin{cases}1&\text{if }x\leq y,\\y&\text{otherwise},\end{cases}$
\item $v(\bigwedge\Phi)=\inf\{v(\phi)\mid \phi\in\Phi\}$,
\item $v(\bigvee\Phi)=\sup\{v(\phi)\mid \phi\in\Phi\}$,
\end{enumerate}
for any $\Phi\cup\{\phi,\psi\}\subseteq\mathcal{L}_\kappa$.
\end{definition}
Given a set of formulas $\Gamma$, we write $v[\Gamma]:=\{v(\gamma)\mid\gamma\in\Gamma\}$. The derived notion of semantics consequence is then defined as follows: for $\Gamma\cup\{\phi\}\subseteq\mathcal{L}_\kappa$, we write $\Gamma\models_\mathsf{G_\kappa}\phi$ if $v[\Gamma]\subseteq\{1\}$ implies $v(\phi)=1$ for any $\mathcal{L}_\kappa$-G\"odel-evaluation $v$.\\

We call the set of consequences $\Gamma\models_\mathsf{G_\kappa}\phi$ the \emph{$\kappa$-infinitary G\"odel logic} and denote it by $\mathsf{G}_\kappa$. 
For the particular case of $\mathsf{G}_{\omega_1}$, a main part of the paper is devoted to the study of various proof theoretic formalism for capturing that semantic consequence and we thus continue by introducing the relevant Hilbert-style calculus used later in a corresponding completeness proof.
\subsection{A Proof Calculus for $\kappa=\omega_1$}
The proof calculus for $\mathsf{G}_{\omega_1}$ which we introduce, denoted by $\mathcal{G}_{\omega_1}$, is a straightforward combination of a proof calculus for propositional infinitary intuitionistic logic with the (pre-)linearity scheme
\[
(\phi\rightarrow\psi)\lor (\psi\rightarrow\phi).
\]
To be concrete, we consider the following system of axiom and rules:
\subsection*{The Calculus $\mathcal{G}_{\omega_1}$}
\begin{description}
\item [($IL$)] a complete set of axioms for propositional intuitionistic logic;\footnote{Naturally, we expect the set of axioms to be defined using the connectives $\land,\lor,\rightarrow,\bot$ and to only require modus ponens as an inference rule. For a particular choice, take the schemes $\land$-Ax, $\lor$-Ax, $\rightarrow$-Ax, $\bot$-Ax from \cite{TvD1988}.}
\item [($GL$)] $(\phi\rightarrow\psi)\lor (\psi\rightarrow\phi)$;
\item [($\omega\lor$)] $\phi_j\rightarrow\bigvee_{i\in\omega}\phi_i$, ($j\in\omega$);
\item [($\omega\land$)] $\bigwedge_{i\in\omega}\phi_i\rightarrow\phi_j$, ($j\in\omega$);
\item [($MP$)] from $\phi\rightarrow\psi$ and $\phi$, infer $\psi$;
\item [($R\omega$)$_1$] from $\phi_i\rightarrow\psi$ for all $i\in\omega$, infer $\bigvee_{i\in\omega}\phi_i\rightarrow\psi$;
\item [($R\omega$)$_2$] from $\phi\rightarrow\psi_i$ for all $i\in\omega$, infer $\phi\rightarrow\bigwedge_{i\in\omega}\psi_i$.
\end{description}
Further, we consider the extension $\mathcal{G}^D_{\omega_1}$ which extends the above calculus by the axiom scheme
\[
\tag{D}\bigwedge_{i\in\omega}(\phi\lor\psi_i)\rightarrow\left(\phi\lor\bigwedge_{i\in\omega}\psi_i\right)
\]
expressing the distributivity of the infinitary operations.
A proof in $\mathcal{G}_{\omega_1}$ (or $\mathcal{G}^D_{\omega_1}$) of some $\phi$ from some assumptions $\Gamma$ is any function $f:\alpha+1\to\mathcal{L}_{\omega_1}$ where $\alpha<\omega_1$ as well as $f(\alpha)=\phi$ and such that any $f(\beta)$ is either
\begin{enumerate}
\item an instance of an axiom scheme,
\item element of $\Gamma$, 
\item the result of ($MP$) or ($R\omega$)$_1$, ($R\omega$)$_2$ with assumptions $f(\gamma)$ where $\gamma<\beta$.
\end{enumerate}
We write $\Gamma\vdash_{\mathcal{G}_{\omega_1}}\phi$ (or $\Gamma\vdash_{\mathcal{G}^D_{\omega_1}}\phi$, respectively) if there is such a proof.

Relative to some fragment $\mathcal{L}_A$ as defined above, we also introduce a restricted notion of derivation: for $\Gamma\cup\{\phi\}\subseteq\mathcal{L}_A$, we write $\Gamma\vdash_{\mathcal{G}_{\omega_1}(\mathcal{L}_A)}\phi$ (or $\Gamma\vdash_{\mathcal{G}^D_{\omega_1}(\mathcal{L}_A)}\phi$, respectively) if there is a proof $f$ with $\mathrm{img}(f)\subseteq\mathcal{L}_A$.\\

Note that $\mathcal{G}_{\omega_1}(\mathcal{L}_A)$, and thus also $\mathcal{G}^D_{\omega_1}(\mathcal{L}_A)$, have the classical Deduction Theorem.
\begin{lemma}
For any $\Gamma\cup\{\phi,\psi\}\subseteq\mathcal{L}_A$, we have 
\[
\Gamma\cup\{\phi\}\vdash_{\mathcal{G}_{\omega_1}(\mathcal{L}_A)}\psi\text{ iff }\Gamma\vdash_{\mathcal{G}_{\omega_1}(\mathcal{L}_A)}\phi\rightarrow\psi.
\]
The same holds for $\mathcal{G}^D_{\omega_1}(\mathcal{L}_A)$.
\end{lemma}
\section{First-Order Infinitary G\"odel Logics}
For the first-order variant, we assume a standard underlying first-order signature $\sigma$ consisting of any number of predicate symbols $P$ and functions symbols $f$. For any such given symbol, we write $\mathrm{ar}(P)$ or $\mathrm{ar}(f)$, respectively, for its arity (which is assumed to be finite). 

We construct a infinitary language corresponding to cardinals $\kappa\geq\lambda$ as it is usually done classically as well: assume a set of variables of size $\kappa$, given by $Var_{\kappa}$ as before, and in that context denote the set of terms over $\sigma$ and $Var_\kappa$ by $\mathcal{T}_\kappa(\sigma)$. The \emph{infinitary language over $\sigma$ associated with $\kappa,\lambda$} is then given by
\[
\mathcal{L}_{\kappa,\lambda}(\sigma):\phi::=\bot\mid P(t_1,\dots,t_n)\mid (\phi\rightarrow\phi)\mid (\phi\land\phi)\mid(\phi\lor\phi)\mid\bigwedge\Phi\mid\bigvee\Phi\mid\exists X\phi\mid\forall X\phi
\]
where $X\subseteq Var_{\kappa}$ is a set of size $<\lambda$, $P\in\sigma$ is a predicate symbol with $\mathrm{ar}(P)=n$, $t_1,\dots,t_n\in\mathcal{T}_\kappa(\sigma)$ and $\Phi$ is a set of formulas of size $<\kappa$. We drop the $\sigma$ if the context is clear or if the choice is arbitrary.

We write $\mathrm{free}(\phi)$ for the set of free variables of $\phi$ and $\mathrm{var}(\phi)$ for the set of variables of $\phi$, free or bound. As before, we write $\mathrm{sub}(\phi)$ for the set of subformulas of $\phi$. Further, we write $\sigma(\phi)$ for the set of function and predicate symbols occurring in $\phi$. These notions straightforwardly extend to sets $\Gamma$ and we use the same notation there.

Again, the countable case $\mathcal{L}_{\omega_1,\omega}$ with finite quantifiers will be of particular interest, especially in the context of a completeness theorem, later on. In that case, we consider the existential and universal quantifiers to just quantify one variable at a time and write
\[
\exists x\phi\text{ or }\forall x\phi
\]
for $x\in Var_{\omega_1}$ as usual in that case. Further, we again write $\bigwedge_{i\in\omega}\phi_i$ for $\bigwedge\{\phi_i\mid i\in\omega\}$ and $\bigvee_{i\in\omega}\phi_i$ for $\bigvee\{\phi_i\mid i\in\omega\}$ where $(\phi_i)_{i\in\omega}$ is a countable family of formulas. 

On $\mathcal{L}_{\kappa,\lambda}$, we denote the simultaneous substitution of terms $\overline{t}=(t_1,\dots,t_n)$ for free variables $\overline{x}=(x_{i_1},\dots,x_{i_n})$ with $i_j\neq i_k$ for $j\neq k$ in a term $t$ by $t[\overline{t}/\overline{x}]$ and by $\phi[\overline{t}/\overline{x}]$ for formulas $\phi$. Here, we assume that quantifiers are treated by renaming the quantified variable in the sense of
\[
(Qx\phi)[\overline{t}/\overline{x}]:=Qz\phi[(\overline{t}',z)/(\overline{x}',x)]
\]
where $Q\in\{\forall,\exists\}$ and $\overline{x}'$ is $\overline{x}$ with $x$ removed (if it occurs), $\overline{t}'$ is $\overline{t}$ with $t_j$ removed when $x_{i_j}=x$ and $z$ is fresh, i.e. does not occur in $\phi$ or $\overline{t}$.

Similarly to the propositional case, we consider a notion of fragments for $\mathcal{L}_{\omega_1,\omega}$ by extending the previous properties appropriately to allow for a Lindenbaum-Tarski construction over these fragments also in the first-order case later on.
\begin{definition}
A \emph{(distributive) fragment of $\mathcal{L}_{\omega_1,\omega}$} is a set $\mathcal{L}_A\subseteq\mathcal{L}_{\omega_1,\omega}$ together with a set $Var_A\subseteq Var_{\omega_1}$ and a signature $\sigma_A$ such that $\mathcal{T}_A$ is the set of terms of $\sigma_A$ using $Var_A$ and
\begin{enumerate}
\item $\bot\in\mathcal{L}_A$,
\item $P(t_1,\dots, t_n)\in\mathcal{L}_A$ for n-ary $P\in\sigma_A$ and $t_i\in\mathcal{T}_A$,
\item $\phi\in\mathcal{L}_A$ implies $\mathrm{sub}(\phi)\subseteq\mathcal{L}_A$, $\mathrm{var}(\phi)\subseteq Var_A$ and $\sigma(\phi)\subseteq\sigma_A$,
\item $\phi,\psi\in\mathcal{L}_A$ implies $\phi\circ\psi\in\mathcal{L}_A$ for $\circ\in\{\land,\lor,\rightarrow\}$ and $\exists x\phi,\forall x\phi\in\mathcal{L}_A$ for $x\in Var_A$,
\item $\phi\in\mathcal{L}_A$, $t\in\mathcal{T}_A$ implies $\phi[\overline{t}/\overline{x}]\in\mathcal{L}_A$, $t[\overline{t}/\overline{x}]\in\mathcal{T}_A$ for any $\overline{t}\in(\mathcal{T}_A)^n$, any $\overline{x}=(x_{i_1},\dots,x_{i_n})\in (Var_A)^n$ with $i_j\neq i_k$ for $j\neq k$.
\item $\phi,\bigwedge_{i\in\omega}\phi_i\in\mathcal{L}_A$ implies $\bigwedge_{i\in\omega}(\phi\rightarrow\phi_i)\in\mathcal{L}_A$,
\item $\phi,\bigvee_{i\in\omega}\phi_i\in\mathcal{L}_A$ implies $\bigwedge_{i\in\omega}(\phi_i\rightarrow\phi)\in\mathcal{L}_A$,
\item $\phi,\bigwedge_{i\in\omega}\phi_i\in\mathcal{L}_A$ implies $\bigwedge_{i\in\omega}(\phi\lor\phi_i)\in\mathcal{L}_A$.
\end{enumerate}
It is additionally assumed that fragments are ``saturated" when it comes to variables, in the sense that there are enough variables to find fresh ones given a finite selection of formulas from $\mathcal{L}_A$. More precisely, we want that for any $\phi_1,\dots,\phi_n\in\mathcal{L}_A$, there is a variable $y\not\in\mathrm{var}(\phi_1)\cup\dots\cup\mathrm{var}(\phi_n)$.
\end{definition}
\begin{lemma}\label{lem:genfrag}
For any $\Gamma\subseteq\mathcal{L}_{\omega_1,\omega}$, there is a smallest (w.r.t. $\subseteq$) distributive fragment $\mathrm{frag}(\Gamma)$ such that $\Gamma\subseteq\mathrm{frag}(\Gamma)$. If $\Gamma$ is countable, then $\mathrm{frag}(\Gamma)$ is also countable.
\end{lemma}
Note that, for countable $\Gamma$, one can even find a countable fragment $\mathcal{L}_A\supseteq\Gamma$ with a countable $Y\subseteq Var_A$ such that $\mathrm{var}(\phi)\cap Y$ is finite for any $\phi\in\mathcal{L}_A$. So in this case, the saturation of variables is directly satisfied.
\subsection{The Standard Semantics and $\mathsf{G}_{\kappa,\lambda}$}
The standard semantics of infinitary first-order G\"odel logics which we want to consider is, like in the propositional case, a straightforward extension of the usual finitary case. Also here, we initially focus on the full real unit interval $[0,1]$ as the corresponding truth-value set.
\begin{definition}
An \emph{$\mathcal{L}_{\kappa,\lambda}(\sigma)$-G\"odel-model} is a structure $\mathfrak{M}$ which consists of 
\begin{enumerate}
\item a non-empty set $M$,
\item $P^\mathfrak{M}:M^n\to [0,1]$ for every $n$-ary predicate $P$ of $\sigma$,
\item $f^\mathfrak{M}:M^n\to M$ for every $n$-ary function $f$ of $\sigma$.
\end{enumerate}
An \emph{$\mathcal{L}_{\kappa,\lambda}(\sigma)$-G\"odel-interpretation} is a structure $\mathfrak{I}=(\mathfrak{M},v)$ composed of an $\mathcal{L}_{\kappa,\lambda}(\sigma)$-G\"odel-model $\mathfrak{M}$ together with a function $v:Var_\kappa\to M$.
\end{definition}
Over such an interpretation $\mathfrak{I}$, one naturally defines the value $t^\mathfrak{I}$ of some term $t$ of $\mathcal{T}_\kappa$. Further, we define
\[
v\tfrac{f}{X}(x):=\begin{cases}f(x)&\text{if }x\in X,\\v(x)&\text{otherwise},\end{cases}
\]
for $X\subseteq Var_\kappa$ and functions $f:X\to M$. We write $v\tfrac{m}{x}$ for the case of $X=\{x\}$ and $f(x)=m$ and also introduce a special notation for finite tuples with
\[
v\tfrac{\overline{m}}{\overline{x}}:=\left(\dots\left(v\tfrac{m_1}{x_1}\right)\dots\right)\tfrac{m_n}{x_n}
\]
where $\overline{m}=(m_1,\dots,m_n)\in M^n$ and $\overline{x}=(x_{i_1},\dots,x_{i_n})\in (Var_\kappa)^n$. We write
\[
\mathfrak{I}\tfrac{f}{X}:=(\mathfrak{M},v\tfrac{f}{X})
\]
and similarly for singletons and tuples. We also allow empty sets/tuples $\overline{m}$, $\overline{x}$ and set $v\tfrac{\overline{m}}{\overline{x}}:=v$ in this case. 

By recursion on $\mathcal{L}_{\kappa,\lambda}$, we construct the \emph{evaluation $\overline{\mathfrak{I}}:\mathcal{L}_{\kappa,\lambda}\to [0,1]$ associated with} $\mathfrak{I}$:
\begin{enumerate}
\item $\overline{\mathfrak{I}}(\bot):=0^\mathbf{A}$;
\item $\overline{\mathfrak{I}}(P(t_1,\dots,t_n)):=P^\mathfrak{M}(t_1^\mathfrak{I},\dots,t_n^\mathfrak{I})$ for $n$-ary $P$;
\item $\overline{\mathfrak{I}}(\phi\land\psi):=\min\{\overline{\mathfrak{I}}(\phi),\overline{\mathfrak{I}}(\psi)\}$;
\item $\overline{\mathfrak{I}}(\phi\lor\psi):=\max\{\overline{\mathfrak{I}}(\phi),\overline{\mathfrak{I}}(\psi)\}$;
\item $\overline{\mathfrak{I}}(\phi\rightarrow\psi):=\overline{\mathfrak{I}}(\phi)\Rightarrow \overline{\mathfrak{I}}(\psi)$;
\item $\overline{\mathfrak{I}}\left(\bigwedge\Phi\right):=\inf\{\overline{\mathfrak{I}}(\phi)\mid \phi\in\Phi\}$;
\item $\overline{\mathfrak{I}}\left(\bigvee\Phi\right):=\sup\{\overline{\mathfrak{I}}(\phi)\mid\phi\in\Phi\}$;
\item $\overline{\mathfrak{I}}(\forall X\phi):=\inf\{\overline{\mathfrak{I}\tfrac{f}{X}}(\phi)\mid f:X\to M\}$;
\item $\overline{\mathfrak{I}}(\exists X\phi):=\sup\{\overline{\mathfrak{I}\tfrac{f}{X}}(\phi)\mid f:X\to M\}$.
\end{enumerate}
As before, one immediately derives a notion of semantical consequence from the model/interpretation construction and their corresponding evaluations: for $\Gamma\cup\{\phi\}\subseteq\mathcal{L}_{\kappa,\lambda}$, we write $\Gamma\models_\mathsf{G_{\kappa,\lambda}}\phi$ if $\overline{\mathfrak{I}}[\Gamma]\subseteq\{1\}$ implies $\overline{\mathfrak{I}}(\phi)=1$ for any $\mathcal{L}_{\kappa,\lambda}$-G\"odel-interpretation $\mathfrak{I}$.\\

We again define the $\kappa,\lambda$-infinitary G\"odel logic to be the set of consequences $\Gamma\models_\mathsf{G_{\kappa,\lambda}}\phi$ and in general denote it by $\mathsf{G_{\kappa,\lambda}}$.
\subsection{A Proof Calculus}
The following proof calculus $\mathcal{G}_{\omega_1,\omega}$ is obtained by extending the previous proof calculus for the propositional case with appropriate axioms and rules for the quantifiers.
\subsection*{The Calculus $\mathcal{G}_{\omega_1,\omega}$}
\begin{description}
\item [($IL$)] a complete set of axiom schemes for propositional intuitionistic logic, in the first-order language;\footnote{The same remark as in Section 2.3 applies, in fact one can again just take the schemes $\land$-Ax, $\lor$-Ax, $\rightarrow$-Ax, $\bot$-Ax from \cite{TvD1988}, now in the first-order language.}
\item [($GL$)] $(\phi\rightarrow\psi)\lor (\psi\rightarrow\phi)$;
\item [($\omega\lor$)] $\phi_j\rightarrow\bigvee_{i\in\omega}\phi_i$, ($j\in\omega$);
\item [($\omega\land$)] $\bigwedge_{i\in\omega}\phi_i\rightarrow\phi_j$, ($j\in\omega$);
\item [($\forall E$)] $\forall x \phi\rightarrow \phi[t/x]$;
\item [($\exists E$)] $\phi[t/x]\rightarrow\exists x\phi$;
\item [($MP$)] from $\phi\rightarrow\psi$ and $\phi$, infer $\psi$;
\item [($R\omega$)$_1$] from $\phi_i\rightarrow\psi$ for $i\in\omega$, infer $\bigvee_{i\in\omega}\phi_i\rightarrow\psi$;
\item [($R\omega$)$_2$] from $\phi\rightarrow\psi_i$ for $i\in\omega$, infer $\phi\rightarrow\bigwedge_{i\in\omega}\psi_i$;
\item [($\forall I$)] from $\psi\rightarrow\phi$, infer $\psi\rightarrow\forall x\phi$ where $x\not\in\mathrm{free}(\psi)$;
\item [($\exists I$)] from $\phi\rightarrow\psi$, infer $\exists x\phi\rightarrow\psi$ where $x\not\in\mathrm{free}(\psi)$.
\end{description}
As before, we consider an extension $\mathcal{G}^D_{\omega_1,\omega}$ obtained by adding the scheme
\[
\tag{D}\bigwedge_{i\in\omega}(\phi\lor\psi_i)\rightarrow\left(\phi\lor\bigwedge_{i\in\omega}\psi_i\right)
\]
and now \emph{additionally} also the axiom scheme
\[
\tag{QS}\forall x(\psi\lor\phi)\rightarrow (\psi\lor\forall x\phi)\text{ where }x\not\in\mathrm{free}(\psi).
\]
The notion of proof immediately transfers to this setting from the propositional case. We write $\Gamma\vdash_{\mathcal{G}_{\omega_1,\omega}}\phi$ (or $\Gamma\vdash_{\mathcal{G}^D_{\omega_1,\omega}}\phi$, respectively) if there is such a proof. We define restrictions $\mathcal{G}_{\omega_1,\omega}(\mathcal{L}_A)$ (or $\mathcal{G}^D_{\omega_1,\omega}(\mathcal{L}_A)$) to some fragment $\mathcal{L}_A$ of $\mathcal{L}_{\omega_1,\omega}$ as before.\\

Note that also both $\mathcal{G}_{\omega_1,\omega}(\mathcal{L}_A)$ and $\mathcal{G}^D_{\omega_1,\omega}(\mathcal{L}_A)$ have the classical Deduction Theorem.
\begin{lemma}
For any $\Gamma\cup\{\phi,\psi\}\subseteq\mathcal{L}_A$ and $\phi$ closed, we have 
\[
\Gamma\cup\{\phi\}\vdash_{\mathcal{G}_{\omega_1,\omega}(\mathcal{L}_A)}\psi\text{ iff }\Gamma\vdash_{\mathcal{G}_{\omega_1,\omega}(\mathcal{L}_A)}\phi\rightarrow\psi.
\]
The same holds for $\mathcal{G}^D_{\omega_1,\omega}(\mathcal{L}_A)$.
\end{lemma}
\section{L-Algebras, Chains and Algebraic Semantics}
We follow a similar route to semantic completeness as in the setting of finitary G\"odel logics (see in particular \cite{BPZ2007,Hor1969}): we first establish completeness w.r.t. a class of algebras and then construct embeddings from that class into the relevant structures of the intended interpretation. 

More precisely, we first show completeness w.r.t. linearly ordered and sufficiently complete Heyting algebras over countable fragments and then extend this to the Heyting algebra of the real unit interval by embeddings, similar to \cite{Hor1969}. This approach does not only offer a high degree of modularity but also establishes linear Heyting algebras (with sufficient completeness) as the algebraic semantics for infinitary G\"odel logics, in analogy to the finitary case. Once we have established the result with respect to countable fragments, this assumption can be removed over complete algebras like the unit interval.

For that, we need various notions from the theory of Heyting algebras and the next subsection gives a, for reasons of self-containedness, quite detailed account mostly following \cite{Hor1969,RS1963} (up to some notation change).
\subsection{Heyting algebras and related notions}
A \emph{Heyting algebra} is a structure $\mathbf{A}$ = $\langle A,\land^\mathbf{A},\lor^\mathbf{A},\rightarrow^\mathbf{A},0^\mathbf{A},1^\mathbf{A}\rangle$ such that $\langle A,\land^\mathbf{A},\lor^\mathbf{A},0^\mathbf{A},1^\mathbf{A}\rangle$ is a bounded lattice with largest element $1^\mathbf{A}$ and smallest element $0^\mathbf{A}$ and $\rightarrow^\mathbf{A}$ is a binary operation with
\begin{enumerate}
\item $x\rightarrow^\mathbf{A}x=1^\mathbf{A}$,
\item $x\land^\mathbf{A}(x\rightarrow^\mathbf{A} y)=x\land^\mathbf{A} y$,
\item $y\land^\mathbf{A}(x\rightarrow^\mathbf{A} y)=y$,
\item $x\rightarrow^\mathbf{A} (y\land^\mathbf{A} z)=(x\rightarrow^\mathbf{A} y)\land^\mathbf{A}(x\rightarrow^\mathbf{A} z)$,
\end{enumerate}
where we write $a\leq^\mathbf{A}b$ for $a\land^\mathbf{A}b=a$ and $\neg^\mathbf{A}x:=x\rightarrow^\mathbf{A}0^\mathbf{A}$. \emph{Joins} (\emph{suprema}) and \emph{meets} (\emph{infima}) of subsets $X$ are defined as usual and denoted $\Vee{\mathbf{A}}{} X$ and $\Wedge{\mathbf{A}}{} X$, respectively. If every subset has a join and meet, $\mathbf{A}$ is called \emph{complete}. An existing meet $\Wedge{\mathbf{A}}{}X$ is called \emph{distributive} if
\[
\Wedge{\mathbf{A}}{x\in X}(y\lor^\mathbf{A}x)=y\lor^\mathbf{A}\Wedge{\mathbf{A}}{}X
\]
for any $y\in\mathbf{A}$ and $\mathbf{A}$ is called \emph{distributive} if every meet is distributive.

Two particular types of Heyting algebras, which we will consider in this note are \emph{chains}, i.e. Heyting algebras where $\leq^\mathbf{A}$ is linear, and \emph{L-algebras}, i.e. Heyting algebras where $(x\rightarrow^\mathbf{A}y)\lor^\mathbf{A}(y\rightarrow^\mathbf{A}x)=1^\mathbf{A}$ for all $x,y\in\mathbf{A}$. We denote the class of all L-algebras by $\mathsf{L}$, the class of all distributive L-algebras by $\mathsf{DL}$ and the class of all chains by $\mathsf{C}$. We write $\mathsf{CC}$ for the class of countable chains. Naturally, every chain is a distributive L-algebra.\\

Further, we will need the notion of a filter. A set $F\subseteq\mathbf{A}$ is a \emph{filter} for a Heyting algebra $\mathbf{A}$, if (1) $1\in F$, (2) $x,y\in F$ implies $x\land^\mathbf{A}y\in F$ and (3) $x\in F$ and $x\leq^\mathbf{A}y$ imply $y\in F$. $F$ is called \emph{proper} if $F\subsetneq\mathbf{A}$ and $F$ is called a \emph{prime filter} if it is proper and if $x\lor^\mathbf{A}y\in F$ implies $x\in F$ or $y\in F$.\\

We then can ``filter" a Heyting algebra $\mathbf{A}$ via $F$: define $x\leq_F y$ if $x\rightarrow^\mathbf{A}y\in F$ for $x,y\in\mathbf{A}$ and $x\equiv_F y$ if $x\leq_F y$ and $y\leq_Fx$. Then $\equiv_F$ is a congruence relation on Heyting algebras and thus defines a quotient Heyting algebra $\mathbf{A}/F$ over the set of equivalence classes $[a]_F$ of elements $a$ of $\mathbf{A}$ over $\equiv_F$.
In particular, ``filtering" with a prime filter in L-algebras yields a chain:
\begin{lemma}
If $F$ is a prime filter of some L-algebra $\mathbf{A}$, then $\mathbf{A}/F$ is a chain. 
\end{lemma}
The behavior of meets and joins under quotients will be of particular importance later on. For that, we first note the following:
\begin{lemma}[\cite{RS1963}, Chapter IV, 7.2, (7) and (8)]
Let $\mathbf{A}$ be a Heyting algebra. If
\[
\Wedge{\mathbf{A}}{}X\text{ and }\Vee{\mathbf{A}}{}Y
\]
exist in $\mathbf{A}$, then
\[
\Wedge{\mathbf{A}}{x\in X}(z\rightarrow^\mathbf{A}x)=z\rightarrow^\mathbf{A}\Wedge{\mathbf{A}}{}X\text{ and }\Wedge{\mathbf{A}}{y\in Y}(y\rightarrow^\mathbf{A} z)=\Vee{\mathbf{A}}{}Y\rightarrow^\mathbf{A}z
\]
for any $z\in\mathbf{A}$.
\end{lemma}
A filter $F$ of a Heyting algebra $\mathbf{A}$ is said to \emph{preserve an existing meet}
\[
\overline x=\Wedge{\mathbf{A}}{}X
\]
if $\overline x\in F$ if, and only if $x\in F$ for all $x\in X$.

Further, a homomorphism $h:\mathbf{A}\to\mathbf{B}$ of Heyting algebras is said to \emph{preserve a meet} $\Wedge{\mathbf{A}}{}X$, \emph{or a join} $\Vee{\mathbf{A}}{}Y$, if
\[
h\left(\Wedge{\mathbf{A}}{}X\right)=\Wedge{\mathbf{B}}{}h[X]\text{ or }h\left(\Vee{\mathbf{A}}{}Y\right)=\Vee{\mathbf{B}}{}h[Y],
\]
respectively.
\begin{lemma}[\cite{Hor1969}, Lemma 2.2]\label{lem:meetjoinfilter}
Suppose 
\[
\Wedge{\mathbf{A}}{}X\text{ and }\Vee{\mathbf{A}}{}Y
\]
exist in $\mathbf{A}$ and $F$ is a filter of $\mathbf{A}$ which preserves 
\[
\Wedge{\mathbf{A}}{x\in X}(z\rightarrow^\mathbf{A}x)\text{ and }\Wedge{\mathbf{A}}{y\in Y}(y\rightarrow^\mathbf{A} z)
\]
for any $z\in\mathbf{A}$. Then
\[
\left[\Wedge{\mathbf{A}}{}X\right]_F=\Wedge{\mathbf{A}/F}{x\in X}[x]_F\text{ and }\left[\Vee{\mathbf{A}}{}Y\right]_F=\Vee{\mathbf{A}/F}{y\in Y}[y]_F.
\]
Therefore, the canonical map $x\mapsto [x]_F$ from $\mathbf{A}$ into $\mathbf{A}/F$ is a Heyting algebra homomorphism which preserves the respective meet and join.
\end{lemma}
\begin{lemma}[\cite{Hor1969}, Theorem 2.3]\label{lem:filterlemma}
Let $\mathbf{A}$ be a Heyting algebra and let
\[
x_n=\Wedge{\mathbf{A}}{}X_n
\]
be a sequence of distributive meets in $\mathbf{A}$. If $x,y\in\mathbf{A}$ with $x\not\leq^\mathbf{A}y$, then there is a prime filter $F$ such that $x\in F$, $y\not\in F$ and such that $F$ preserves the meets given by $x_n$.
\end{lemma}
By $\mathbf{[0,1]_\mathbb{Q}}$ and $\mathbf{[0,1]_\mathbb{R}}$, we denote the Heyting algebras of all rationals in the unit interval and of the whole unit interval, respectively.
\begin{lemma}[\cite{Hor1969}, Lemma 3.7]\label{lem:qembed}
Let $\mathbf{A}$ be a countable chain. Then there is an embedding, i.e. an injective homomorphism of Heyting algebras
\[
q:\mathbf{A}\to\mathbf{[0,1]_\mathbb{Q}}
\]
which preserves all meets and joins of $\mathbf{A}$.
\end{lemma}
\subsection{Algebraic Propositional Evaluations for $\mathcal{L}_{\omega_1}$}
In this section, we now introduce the actual algebraic generalizations of the $\mathcal{L}_{\omega_1}$-G\"odel-evaluations, broadening the domains to fragments and the range to certain Heyting algebras which may be, in a particular way, incomplete. This will be necessary in the approach to completeness chosen here since the Lindenbaum-Tarski algebras later constructed are, in fact, incomplete. 

Let $\mathcal{L}_A$ be an arbitrary fragment of $\mathcal{L}_{\omega_1}$ and $\mathbf{A}$ be a Heyting algebra.
\begin{definition}
A function $v:\mathcal{L}_A\to\mathbf{A}$ is an \emph{($\mathbf{A}$-valued) $\mathcal{L}_A$-evaluation} if
\begin{enumerate}
\item $v(\bot)=0^\mathbf{A}$,
\item $v(\phi\circ\psi)=v(\phi)\circ^\mathbf{A}v(\psi)$ for $\circ\in\{\rightarrow,\land,\lor\}$,
\item for any $\bigwedge_{i\in\omega}\phi_i,\bigvee_{i\in\omega}\psi_i\in\mathcal{L}_A$, we have
\[
v\left(\bigwedge_{i\in\omega}\phi_i\right)=\Wedge{\mathbf{A}}{i\in\omega}v(\phi_i)
\]
and 
\[
v\left(\bigvee_{i\in\omega}\phi_i\right)=\Vee{\mathbf{A}}{i\in\omega}v(\phi_i)
\]
such that the corresponding infima/suprema exist.
\end{enumerate}
\end{definition}
Given such an evaluation $v$, we still write $v[\Gamma]:=\{v(\gamma)\mid\gamma\in\Gamma\}$ for sets $\Gamma\subseteq\mathcal{L}_{A}$ and $(\mathbf{A},v)\models\phi$ for $v(\phi)=1^\mathbf{A}$. We denote the set of all $\mathbf{A}$-valued $\mathcal{L}_A$-evaluations by $\mathsf{Ev}(\mathcal{L}_A;\mathbf{A})$.\\

Using this notion of $\mathcal{L}_A$-evaluations, there is now a natural notion of semantical entailment: let $\mathsf{Cl}$ be a class of Heyting algebras and let $\Gamma\cup\{\phi\}\subseteq\mathcal{L}_A$. We write $\Gamma\models_{\mathsf{Cl}(\mathcal{L}_A)}\phi$ if
\[
\forall\mathbf{A}\in\mathsf{Cl}\forall v\in\mathsf{Ev}(\mathcal{L}_A;\mathbf{A})\left(v[\Gamma]\subseteq\{1^\mathbf{A}\}\text{ implies }v(\phi)=1^\mathbf{A}\right).
\]
We abbreviate $\Gamma\models_{\mathsf{Cl}(\mathcal{L}_{\omega_1})}\phi$ by $\Gamma\models_{\mathsf{Cl}}\phi$.\\

By transfinite induction on the length of the proof, on quickly verifies the following soundness result:
\begin{lemma}\label{lem:sound}
For any $\Gamma\cup\{\phi\}\subseteq\mathcal{L}_A$:
\begin{enumerate}
\item $\Gamma\vdash_{\mathcal{G}_{\omega_1}(\mathcal{L}_A)}\phi$ implies $\Gamma\models_{\mathsf{L}(\mathcal{L}_A)}\phi$,
\item $\Gamma\vdash_{\mathcal{G}^D_{\omega_1}(\mathcal{L}_A)}\phi$ implies $\Gamma\models_{\mathsf{DL}(\mathcal{L}_A)}\phi$.
\end{enumerate}
\end{lemma}
To form new $\mathcal{L}_A$-evaluations by composition with homomorphisms of Heyting algebras, we now have to additionally require that the existing meets and joins arising from the infinitary connectives are preserved. This is captured in the following lemma.
\begin{lemma}\label{lem:homprop}
Let $v:\mathcal{L}_A\to\mathbf{A}$ be an evaluation of $\mathcal{L}_A$ in $\mathbf{A}$ and let $h:\mathbf{A}\to\mathbf{B}$ be a Heyting algebra homomorphism which preserves all the meets and joins
\[
\Wedge{\mathbf{A}}{i\in\omega}v(\phi_i)\text{ and }\Vee{\mathbf{A}}{i\in\omega}v(\psi_i)
\]
for any $\bigwedge_{i\in\omega}\phi_i\in\mathcal{L}_A$ and $\bigvee_{i\in\omega}\psi_i\in\mathcal{L}_A$. Then $h\circ v$ is an evaluation of $\mathcal{L}_A$ in $\mathbf{B}$.
\end{lemma}
The proof is rather immediate and thus omitted.
\subsection{Algebraic First-Order Interpretations for $\mathcal{L}_{\omega_1,\omega}$}
For a similar motivation as in the propositional case, we are also lead to broadening the definitions of models and interpretations for $\mathcal{L}_{\omega_1,\omega}$ to both arbitrary fragments as domains and Heyting algebras as ranges. For this, let $\mathcal{L}_A$ now be a fragment of $\mathcal{L}_{\omega_1,\omega}$ and $\mathbf{A}$ again be a Heyting algebra.
\begin{definition}
An \emph{$\mathcal{L}_A$-model} is a structure $\mathfrak{M}$ which consists of 
\begin{enumerate}
\item a Heyting algebra $\mathbf{A}$,
\item a non-empty set $M$,
\item $P^\mathfrak{M}:M^n\to\mathbf{A}$ for every $n$-ary predicate $P$ of $\sigma_A$,
\item $f^\mathfrak{M}:M^n\to M$ for every $n$-ary function $f$ of $\sigma_A$.
\end{enumerate}
An \emph{$\mathcal{L}_A$-interpretation} is a structure $\mathfrak{I}=(\mathfrak{M},v)$ where $\mathfrak{M}$ is an $\mathcal{L}_A$-model and $v:Var_A\to M$.
\end{definition}
As before, over some $\mathcal{L}_A$-interpretation $\mathfrak{I}$, one naturally defines the value of some term $t$ of $\mathcal{T}_A$ which we still denote by $t^\mathfrak{I}$. $v\tfrac{m}{x}$, $v\tfrac{\overline{m}}{\overline{x}}$ and the resulting $\mathfrak{I}\tfrac{\overline{m}}{\overline{x}}$ are defined (over $Var_A$) in the same way as with the standard semantics.

Such a model $\mathfrak{M}$ is called \emph{($\mathcal{L}_A$-)suitable w.r.t. to a variable assignment} $v:Var_A\to M$ if for any $\overline{m}\in M^n$, any $\overline{x}\in (Var_A)^n$, corresponding to the interpretation $\mathfrak{I}=(\mathfrak{M},v\tfrac{\overline{m}}{\overline{x}})$ there is a function $\overline{\mathfrak{I}}:\mathcal{L}_A\to\mathbf{A}$ such that
\begin{enumerate}
\item $\overline{\mathfrak{I}}(\bot)=0^\mathbf{A}$,
\item $\overline{\mathfrak{I}}(P(t_1,\dots,t_n))=P^\mathfrak{M}(t_1^\mathfrak{I},\dots,t_n^\mathfrak{I})$ for $n$-ary $P$,
\item $\overline{\mathfrak{I}}(\phi\circ\psi)=\overline{\mathfrak{I}}(\phi)\circ^\mathbf{A}\overline{\mathfrak{I}}(\psi)$ for $\circ\in\{\land,\lor,\rightarrow\}$,
\item for any $\bigwedge_{i\in\omega}\phi_i,\bigvee_{i\in\omega}\psi_i\in\mathcal{L}_A$, we have
\[
\overline{\mathfrak{I}}\left(\bigwedge_{i\in\omega}\phi_i\right)=\Wedge{\mathbf{A}}{i\in\omega}\overline{\mathfrak{I}}(\phi_i)
\]
and 
\[
\overline{\mathfrak{I}}\left(\bigvee_{i\in\omega}\phi_i\right)=\Vee{\mathbf{A}}{i\in\omega}\overline{\mathfrak{I}}(\phi_i)
\]
such that the corresponding infima/suprema exist,
\item for any $\phi\in\mathcal{L}_A$ and any $x\in Var_A$, we have
\[
\overline{\mathfrak{I}}(\forall x\phi)=\Wedge{\mathbf{A}}{m\in M}\overline{\mathfrak{I}\tfrac{m}{x}}(\phi)
\]
and
\[
\overline{\mathfrak{I}}(\exists x\phi)=\Vee{\mathbf{A}}{m\in M}\overline{\mathfrak{I}\tfrac{m}{x}}(\phi)
\]
such that the corresponding infima/suprema exist.
\end{enumerate}
\begin{remark}\label{rem:suitableshortcut}
Note that it actually suffices to establish the existence of such an extension only for $\overline{m}$ and $\overline{x}=(x_{i_1},\dots,x_{i_n})$ where $i_j\neq i_k$ for $j\neq k$.
\end{remark}
We still write $\overline{\mathfrak{I}}[\Gamma]:=\{\overline{\mathfrak{I}}(\gamma)\mid\gamma\in\Gamma\}$ for sets $\Gamma\subseteq\mathcal{L}_A$ and $\overline{\mathfrak{I}}\models\phi$ if $\overline{\mathfrak{I}}(\phi)=1^\mathbf{A}$. 

We denote the class of all models over $\mathbf{A}$ by $\mathsf{Mod}(\mathcal{L}_A;\mathbf{A})$ and if $\mathfrak{M}$ is a model, we write $\mathsf{Int}(\mathcal{L}_A;\mathfrak{M})$ for the set of all corresponding interpretations $(\mathfrak{M},v)$ such that $\mathfrak{M}$ is suitable for $v$.\\

The derived notion of semantic consequence is then given by the following: let $\mathsf{Cl}$ be a class of Heyting algebras and let $\Gamma\cup\{\phi\}\subseteq\mathcal{L}_A$. We write $\Gamma\models_{\mathsf{Cl}(\mathcal{L}_A)}\phi$ if
\[
\forall\mathbf{A}\in\mathsf{Cl}\forall\mathfrak{M}\in\mathsf{Mod}(\mathcal{L}_A;\mathbf{A})\forall\mathfrak{I}\in\mathsf{Int}(\mathcal{L}_A;\mathfrak{M})\left(\overline{\mathfrak{I}}[\Gamma]\subseteq\{1^\mathbf{A}\}\text{ implies }\overline{\mathfrak{I}}(\phi)=1^\mathbf{A}\right).
\]
It is again straightforward to verify the following soundness result.
\begin{lemma}\label{lem:fosound}
For any $\Gamma\cup\{\phi\}\subseteq\mathcal{L}_A$ where all formulas from $\Gamma$ are closed:
\begin{enumerate}
\item $\Gamma\vdash_{\mathcal{G}_{\omega_1,\omega}(\mathcal{L}_A)}\phi$ implies $\Gamma\models_{\mathsf{L}(\mathcal{L}_A)}\phi$,
\item $\Gamma\vdash_{\mathcal{G}^D_{\omega_1,\omega}(\mathcal{L}_A)}\phi$ implies $\Gamma\models_{\mathsf{DL}(\mathcal{L}_A)}\phi$.
\end{enumerate}
\end{lemma}
As before, we can form new interpretations by composition with Heyting algebra homomorphisms as long as they respect the need meets and joins.
\begin{lemma}\label{lem:fohomprop}
Let $\mathfrak{M}$ be a model over some Heyting algebra $\mathbf{A}$ which is suitable w.r.t. $v:Var_A\to M$ and let $h:\mathbf{A}\to\mathbf{B}$ be a Heyting algebra homomorphism which preserves all the meets and joins
\[
\Wedge{\mathbf{A}}{i\in\omega}\overline{\mathfrak{I}}(\phi_i)\text{ and }\Vee{\mathbf{A}}{i\in\omega}\overline{\mathfrak{I}}(\psi_i)
\]
for any $\bigwedge_{i\in\omega}\phi_i\in\mathcal{L}_A$ and $\bigvee_{i\in\omega}\psi_i\in\mathcal{L}_A$ as well as
\[
\Wedge{\mathbf{A}}{m\in M}\overline{\mathfrak{I}\tfrac{m}{x}}(\phi)\text{ and }\Vee{\mathbf{A}}{m\in M}\overline{\mathfrak{I}\tfrac{m}{x}}(\phi)
\]
for all $\phi\in\mathcal{L}_A$ and $x\in Var_A$ and \emph{for any} interpretation $\mathfrak{I}=(\mathfrak{M},v\tfrac{\overline{m}}{\overline{x}})$ over $\mathfrak{M}$. Then the model $h\circ\mathfrak{M}$ defined over the same domain with $f^{h\circ\mathfrak{M}}:=f^\mathfrak{M}$ for any function symbol $f$ and $P^{h\circ\mathfrak{M}}:=h\circ P^\mathfrak{M}$ for predicate symbols $P$ is a suitable model w.r.t. $v$ and for any interpretation $\mathfrak{I}=(\mathfrak{M},v\tfrac{\overline{m}}{\overline{x}})$, we have
\[
\overline{h\circ\mathfrak{I}}(\phi)=h(\overline{\mathfrak{I}}(\phi))
\]
for all $\phi\in\mathcal{L}_A$ where $h\circ\mathfrak{I}=(h\circ\mathfrak{M},v\tfrac{\overline{m}}{\overline{x}})$.
\end{lemma}
\section{A Propositional and First-Order Completeness Theorem}
We fix a fragment $\mathcal{L}_A$ of either $\mathcal{L}_{\omega_1}$ or $\mathcal{L}_{\omega_1,\omega}$ and all notions, if not explicitly indicated, are to be understood relative to that fragment. Let $\mathcal{G}^{(D)}(\mathcal{L}_A)$ be either $\mathcal{G}_{\omega_1}(\mathcal{L}_A)$ or $\mathcal{G}^D_{\omega_1}(\mathcal{L}_A)$ in the propositional case or either $\mathcal{G}_{\omega_1,\omega}(\mathcal{L}_A)$ or $\mathcal{G}^D_{\omega_1,\omega}(\mathcal{L}_A)$ in the first order case.

We construct, as it is usually done, the Lindenbaum-Tarski algebra of $\mathcal{G}^{(D)}(\mathcal{L}_A)$: given $\Gamma\cup\{\phi,\psi\}\subseteq\mathcal{L}_A$, we write $\phi\preceq^\Gamma\psi$ if
\[
\Gamma\vdash_{\mathcal{G}^{(D)}(\mathcal{L}_A)}\phi\rightarrow\psi
\]
and write $\phi\equiv^\Gamma\psi$ if
\[
\phi\preceq^\Gamma\psi\text{ and }\psi\preceq^\Gamma\phi.
\]
We write $[\phi]^\Gamma$ for the equivalence class of $\phi$ under $\equiv^\Gamma$ and $\mathcal{L}_A/\equiv^\Gamma$ for the set of all equivalence classes. In the following, we even omit the $\Gamma$ when the context is clear. The \emph{Lindenbaum-Tarski algebra} $\mathbf{LT}^\Gamma$ is defined as
\[
\mathbf{LT}^\Gamma:=\langle\mathcal{L}_A/\equiv^\Gamma,\land^{\mathbf{LT}},\lor^\mathbf{LT},\rightarrow^\mathbf{LT},0^\mathbf{LT},1^\mathbf{LT}\rangle
\]
where we define
\begin{enumerate}
\item $[\phi]\land^\mathbf{LT}[\psi]:=[\phi\land\psi]$,
\item $[\phi]\lor^\mathbf{LT}[\psi]:=[\phi\lor\psi]$,
\item $[\phi]\rightarrow^\mathbf{LT}[\psi]:=[\phi\rightarrow\psi]$,
\item $0^\mathbf{LT}:=[\bot]$,
\item $1^\mathbf{LT}:=[\top]$.
\end{enumerate}
Further, the order induced on $\mathbf{LT}^\Gamma$ is given by $[\phi]\leq^\mathbf{LT}[\psi]$ iff $\phi\preceq^\Gamma\psi$
\begin{lemma}\label{lem:ltwelldef}
$\mathbf{LT}^\Gamma$ is a well-defined L-algebra with
\[
\Wedge{\mathbf{LT}}{i\in\omega}[\phi_i]=\left[\bigwedge_{i\in\omega}\phi_i\right]\text{ and }\Vee{\mathbf{LT}}{i\in\omega}[\psi_i]=\left[\bigvee_{i\in\omega}\psi_i\right]
\]
for $\bigwedge_{i\in\omega}\phi_i,\bigvee_{i\in\omega}\psi_i\in\mathcal{L}_A$. In the case of $\mathcal{G}^D(\mathcal{L}_A)$, the meets
\[
[\chi]\rightarrow^\mathbf{LT}\left[\bigwedge_{i\in\omega}\phi_i\right]=\Wedge{\mathbf{LT}}{i\in\omega}([\chi]\rightarrow^\mathbf{LT}[\phi_i])\text{ and }\left[\bigvee_{i\in\omega}\psi_i\right]\rightarrow^\mathbf{LT}[\chi]=\Wedge{\mathbf{LT}}{i\in\omega}([\psi_i]\rightarrow^\mathbf{LT}[\chi])
\]
are distributive for every additional $\chi\in\mathcal{L}_A$. Further, we have
\[
[\forall x\phi]=\Wedge{\mathbf{LT}}{t\in\mathcal{T}_A}\left[\phi[t/x]\right]\text{ and }[\exists x\phi]=\Vee{\mathbf{LT}}{t\in\mathcal{T}_A}\left[\phi[t/x]\right].
\]
for $\phi\in\mathcal{L}_A$ and $x\in Var_A$ in the first-order case and in the case of $\mathcal{G}^D_{\omega_1,\omega}(\mathcal{L}_A)$, the meets
\[
[\chi]\rightarrow^\mathbf{LT}[\forall x\phi]=\Wedge{\mathbf{LT}}{t\in\mathcal{T}_A}([\chi]\rightarrow^\mathbf{LT}\left[\phi[t/x]\right])\text{ and }[\exists x\phi]\rightarrow^\mathbf{LT}[\chi]=\Wedge{\mathbf{LT}}{t\in\mathcal{T}_A}(\left[\phi[t/x]\right]\rightarrow^\mathbf{LT}[\chi]).
\]
are distributive for every additional $\chi\in\mathcal{L}_A$.
\end{lemma}
\begin{proof}
We skip the finitary operations. Using the axiom scheme ($GL$), it is easy to see that $\mathbf{LT}^\Gamma$ is a well-defined L-algebra. We only show the infinitary claims. For that, let $\bigwedge_{i\in\omega}\phi_i\in\mathcal{L}_A$ and $\bigvee_{i\in\omega}\psi_i\in\mathcal{L}_A$. Then, we have
\[
\mathcal{G}^{(D)}(\mathcal{L}_A)\vdash\bigwedge_{i\in\omega}\phi_i\rightarrow\phi_j\text{ and }\mathcal{G}^{(D)}(\mathcal{L}_A)\vdash\psi_j\rightarrow\bigvee_{i\in\omega}\psi_i
\]
for any $j$ by axioms ($\omega\land$) and ($\omega\lor$) which gives 
\[
\left[\bigwedge_{i\in\omega}\phi_i\right]\leq^\mathbf{LT}[\phi_j]\text{ and }[\psi_j]\leq^\mathbf{LT}\left[\bigvee_{i\in\omega}\psi_i\right]
\]
for any $j$. Suppose that $[\chi]\leq^\mathbf{LT}[\phi_j]$ for any $j$ and $[\psi_j]\leq^\mathbf{LT}[\chi']$ for any $j$. Therefore, we have
\[
\Gamma\vdash_{\mathcal{G}^{(D)}(\mathcal{L}_A)}\chi\rightarrow\phi_j\text{ and }\Gamma\vdash_{\mathcal{G}^{(D)}(\mathcal{L}_A)}\psi_j\rightarrow\chi'
\]
for any $j$ which implies
\[
\Gamma\vdash_{\mathcal{G}^{(D)}(\mathcal{L}_A)}\chi\rightarrow\bigwedge_{i\in\omega}\phi_i\text{ and }\Gamma\vdash_{\mathcal{G}^{(D)}(\mathcal{L}_A)}\bigvee_{i\in\omega}\psi_i\rightarrow\chi'
\]
by ($R\omega$)$_{1,2}$ which is
\[
[\chi]\leq^\mathbf{LT}\left[\bigwedge_{i\in\omega}\phi_i\right]\text{ and }
\left[\bigvee_{i\in\omega}\psi_i\right]\leq^\mathbf{LT}[\chi'].
\]
This gives that
\[
\Wedge{\mathbf{LT}}{i\in\omega}[\phi_i]=\left[\bigwedge_{i\in\omega}\phi_i\right]\text{ and }\Vee{\mathbf{LT}}{i\in\omega}[\psi_i]=\left[\bigvee_{i\in\omega}\psi_i\right].
\]

Next, we show that all the mentioned meets are distributive. Let $\chi,\xi\in\mathcal{L}_A$ be arbitrary. We have
\[
\bigwedge_{i\in\omega}(\xi\lor(\psi_i\rightarrow\chi))\in\mathcal{L}_A\text{ and }\bigwedge_{i\in\omega}(\xi\lor(\chi\rightarrow\phi_i))\in\mathcal{L}_A.
\]
by the closure properties of $\mathcal{L}_A$. Write $(\alpha_i)_i$ for either $(\chi\rightarrow\phi_i)_i$ or $(\psi_i\rightarrow\chi)_i$. Now, we get
\[
[\xi]\lor^\mathbf{LT}\Wedge{\mathbf{LT}}{i\in\omega}[\alpha_i]\leq^\mathbf{LT}[\xi]\lor^\mathbf{LT}[\alpha_j]
\]
for any $j$ by axiom $(\omega\land)$ and therefore 
\[
[\xi]\lor^\mathbf{LT}\Wedge{\mathbf{LT}}{i\in\omega}[\alpha_i]\leq^\mathbf{LT}\Wedge{\mathbf{LT}}{i\in\omega}([\xi]\lor^\mathbf{LT}[\alpha_i])=\left[\bigwedge_{i\in\omega}(\xi\lor\alpha_i)\right].
\]
Further, by axiom scheme $(D)$, we have
\[
\Gamma\vdash_{\mathcal{G}^D(\mathcal{L}_A)}\bigwedge_{i\in\omega}(\xi\lor\alpha_i)\rightarrow\left(\xi\lor\bigwedge_{i\in\omega}\alpha_i\right)
\]
which gives the converse
\[
\left[\bigwedge_{i\in\omega}(\xi\lor\alpha_i)\right]\leq^\mathbf{LT}[\xi]\lor^\mathbf{LT}\Wedge{\mathbf{LT}}{i\in\omega}[\alpha_i],
\]
i.e. combined we have
\[
[\xi]\lor^\mathbf{LT}\Wedge{\mathbf{LT}}{i\in\omega}[\alpha_i]=\left[\bigwedge_{i\in\omega}(\xi\lor\alpha_i)\right].
\]

In the first-order case, the quantifier claims can be proved as outlined in \cite{Hor1969}. We still sketch the proof here for self-containedness.

For the first two quantifier claims, note that we have
\[
\Gamma\vdash_{\mathcal{G}^{(D)}_{\omega_1,\omega}(\mathcal{L}_A)}\forall x\phi\rightarrow\phi[t/x]\text{ and }\Gamma\vdash_{\mathcal{G}^{(D)}_{\omega_1,\omega}(\mathcal{L}_A)}\phi[t/x]\rightarrow\exists x\phi
\]
for any $t\in\mathcal{T}_A$ by the axioms $(\forall E)$ and $(\exists E)$. Now, suppose that $\chi\in\mathcal{L}_A$ is such that
\[
\Gamma\vdash_{\mathcal{G}^{(D)}_{\omega_1,\omega}(\mathcal{L}_A)}\chi\rightarrow\phi[t/x]\text{ for all }t\in\mathcal{T}_A\text{ or }\Gamma\vdash_{\mathcal{G}^{(D)}_{\omega_1,\omega}(\mathcal{L}_A)}\phi[t/x]\rightarrow\chi\text{ for all }t\in\mathcal{T}_A.
\]
Then, pick $y\in Var_A$ with $y\not\in\mathrm{var}(\chi)\cup\mathrm{var}(\phi)$. Note that this is possible as fragments are saturated for variables. By assumption, we have
\[
\Gamma\vdash_{\mathcal{G}^{(D)}_{\omega_1,\omega}(\mathcal{L}_A)}\chi\rightarrow\phi[y/x]\text{ or }\Gamma\vdash_{\mathcal{G}^{(D)}_{\omega_1,\omega}(\mathcal{L}_A)}\phi[y/x]\rightarrow\chi
\]
and as $y$ is not free in $\chi$, we get 
\[
\Gamma\vdash_{\mathcal{G}^{(D)}_{\omega_1,\omega}(\mathcal{L}_A)}\chi\rightarrow\forall y\phi[y/x]\text{ or }\Gamma\vdash_{\mathcal{G}^{(D)}_{\omega_1,\omega}(\mathcal{L}_A)}\exists y\phi[y/x]\rightarrow\chi
\]
via the rule $(\forall I)$ and $(\exists I)$. As we have
\[
\vdash_{\mathcal{G}^{(D)}_{\omega_1,\omega}(\mathcal{L}_A)}\forall y\phi[y/x]\leftrightarrow \forall x\phi\text{ and }\vdash_{\mathcal{G}^{(D)}_{\omega_1,\omega}(\mathcal{L}_A)}\exists y\phi[y/x]\leftrightarrow\exists x\phi,
\]
the claims follow.\\

Regarding distributivity, let $\phi\in\mathcal{L}_A$ and $x\in Var_A$ as well as $\chi,\xi\in\mathcal{L}_A$. Further, let $y\in Var_A$ with $y\not\in\mathrm{var}(\phi)\cup\mathrm{var}(\chi)\cup\mathrm{var}(\xi)$ (note again variable saturation). Then, we have
\[
\vdash_{\mathcal{G}^{(D)}_{\omega_1,\omega}(\mathcal{L}_A)}\phi[t/x]\leftrightarrow \phi[y/x][t/y]
\]
and therefore, defining $\phi':=\phi[y/x]$, we have $[\phi[t/x]]=[\phi'[t/y]]$. Since $y$ does not occur in $\chi$, we have
\[
\chi\rightarrow\phi'[t/y]=(\chi\rightarrow\phi')[t/y]=:\alpha^\forall[t/y]
\]
and
\[
\phi'[t/y]\rightarrow\chi=(\phi'\rightarrow\chi)[t/y]=:\alpha^\exists[t/y]
\]
and therefore, since $y$ also does not occur in $\xi$, the previous quantifier claims imply (where we write $\alpha$ for either $\alpha^\exists$ or $\alpha^\forall$):
\begin{align*}
\Wedge{\mathbf{LT}}{t\in\mathcal{T}_A}([\xi]\lor^\mathbf{LT}[\alpha[t/y]])&=\Wedge{\mathbf{LT}}{t\in\mathcal{T}_A}[(\xi\lor\alpha)[t/y]]\\
&=[\forall y(\xi\lor\alpha)].
\end{align*}
Using the distributivity axiom, as $y$ again does not occur in $\xi$, we get 
\begin{align*}
[\forall y(\xi\lor\alpha)]&=[\xi\lor\forall y\alpha]\\
&=[\xi]\lor^\mathbf{LT}\Wedge{\mathbf{LT}}{t\in\mathcal{T}_A}[\alpha[t/y]]
\end{align*}
and finally we have
\[
\Wedge{\mathbf{LT}}{t\in\mathcal{T}_A}[\xi]\lor^\mathbf{LT}([\chi]\rightarrow^\mathbf{LT}\left[\phi[t/x]\right])=[\forall y(\xi\lor\alpha^\forall)]=[\xi]\lor^\mathbf{LT}\Wedge{\mathbf{LT}}{t\in\mathcal{T}_A}([\chi]\rightarrow^\mathbf{LT}\left[\phi[t/x]\right])
\]
and
\[
\Wedge{\mathbf{LT}}{t\in\mathcal{T}_A}[\xi]\lor^\mathbf{LT}(\left[\phi[t/x]\right]\rightarrow^\mathbf{LT}[\chi])=[\forall y(\xi\lor\alpha^\exists)]=[\xi]\lor^\mathbf{LT}\Wedge{\mathbf{LT}}{t\in\mathcal{T}_A}(\left[\phi[t/x]\right]\rightarrow^\mathbf{LT}[\chi]).
\]
\end{proof}
In the propositional case, there is a \emph{canonical $\mathcal{L}_A$-evaluation} over $\mathbf{LT}^\Gamma$, $\iota:\mathcal{L}_A\to\mathbf{LT}^\Gamma$, which is defined by
\[
\iota(\phi):=[\phi].
\]
By Lemma \ref{lem:ltwelldef}, $\iota$ indeed is a well-defined $\mathcal{L}_A$-evaluation.\\

In the first-order case, this Lindenbaum-Tarski algebra now forms the algebraic part of the \emph{Lindenbaum-Tarski model}: define $\mathfrak{LT}^\Gamma$ as a model with $\mathbf{LT}^\Gamma$ as the underlying Heyting algebra and with $\mathcal{T}_A$ as the domain by setting
\[
f^\mathfrak{LT}(t_1,\dots,t_n):=f(t_1,\dots,t_n)\text{ as well as }P^\mathfrak{LT}(t_1,\dots,t_n):=[P(t_1,\dots,t_n)]
\]
for functions symbols $f$ and predicate symbols $P$ of $\mathcal{L}_A$. We define the \emph{canonical variable assignment}
\[
\iota:Var_A\to\mathcal{T}_A,\; x\mapsto x
\]
and denote it, for convenience, also by $\iota$ but the context will make it clear whether we mean the propositional evaluation or the first-order variable assignment.

First note that $\mathfrak{LT}^\Gamma$ is indeed suitable for $\iota$.
\begin{lemma}\label{lem:ltsuitable}
The model $\mathfrak{LT}^\Gamma$ is a suitable model w.r.t. $\iota$. In particular, given some $\overline{t}=(t_1,\dots,t_n)\in(\mathcal{T}_A)^n$ and $\overline{x}=(x_{i_1},\dots,x_{i_n})\in (Var_A)^n$ with $i_j\neq i_k$ for $j\neq k$, it holds that
\[
\overline{(\mathfrak{LT}^\Gamma,\iota\tfrac{\overline{t}}{\overline{x}})}(\phi)=[\phi[\overline{t}/\overline{x}]].
\]
for any $\phi\in\mathcal{L}_A$.
\end{lemma}
\begin{proof}
Consider $(\mathfrak{LT}^\Gamma,\iota\tfrac{\overline{t}}{\overline{x}})$ for $\overline{x}$ and $\overline{t}$ as required. By Remark \ref{rem:suitableshortcut}, showing that $\phi\mapsto [\phi[\overline{t}/\overline{x}]]$ is the extension of $(\mathfrak{LT}^\Gamma,\iota\tfrac{\overline{t}}{\overline{x}})$ for those $\overline{t},\overline{x}$ is enough to establish suitability.

Note first that we have 
\[
s^{\left(\mathfrak{LT}^\Gamma,\iota\tfrac{\overline{t}}{\overline{x}}\right)}=s[\overline{t}/\overline{x}]\in\mathcal{T}_A
\]
for any $s\in\mathcal{T}_A$ as by assumption, $\mathcal{T}_A$ is closed under substitution. The theorem is then proved by induction on $\phi$. By the above result for terms, we get
\[
[P(s_1,\dots,s_n)[\overline{t}/\overline{x}]]=[P(s_1[\overline{t}/\overline{x}],\dots,s_n[\overline{t}/\overline{x}])]=P^{(\mathfrak{LT}^\Gamma,\iota\tfrac{\overline{t}}{\overline{x}})}(s_1^{(\mathfrak{LT}^\Gamma,\iota\tfrac{\overline{t}}{\overline{x}})},\dots,s_n^{(\mathfrak{LT}^\Gamma,\iota\tfrac{\overline{t}}{\overline{x}})}).
\]
Further, we in particular have $[\bot[\overline{t}/\overline{x}]]=[\bot]$ and this provides the result for the atomic cases. The cases for $\phi\circ\psi\in\mathcal{L}_A$ with $\circ\in\{\land,\lor,\rightarrow\}$ follow just by noting that $\cdot[\overline{t}/\overline{x}]$ distributes over $\circ$ and by considering the definitions of $\circ^\mathbf{LT}$.

The same also holds true for $\bigwedge_{i\in\omega}\phi_i\in\mathcal{L}_A$ where we have
\[
\left[\left(\bigwedge_{i\in\omega}\phi_i\right)[\overline{t}/\overline{x}]\right]=\left[\bigwedge_{i\in\omega}\phi_i[\overline{t}/\overline{x}]\right]=\Wedge{\mathbf{LT}}{i\in\omega}\left[\phi_i[\overline{t}/\overline{x}]\right]
\]
and similarly for $\bigvee_{i\in\omega}\phi_i\in\mathcal{L}_A$ where we have used Lemma \ref{lem:ltwelldef}. This gives the infinitary cases.

Lastly, let $\phi\in\mathcal{L}_A$ and $x\in Var_A$. Then, we have
\[
[(\forall x\phi)[\overline{t}/\overline{x}]]=[\forall z\phi[(\overline{t'},z)/(\overline{x'},x)]]
\]
where $z$ and $\overline{x'},\overline{t'}$ are as in the definition of first-order substitutions. Using Lemma \ref{lem:ltwelldef}, we get
\[
[\forall z\phi[(\overline{t'},z)/(\overline{x'},x)]]=\Wedge{\mathbf{LT}}{t\in\mathcal{T}_A}[\phi[(\overline{t'},z)/(\overline{x'},x)][t/z]]=\Wedge{\mathbf{LT}}{t\in\mathcal{T}_A}[\phi[(\overline{t'},t)/(\overline{x'},x)]].
\]
Using similar reasoning, one can show 
\[
[(\exists x\phi)[\overline{t}/\overline{x}]]=\Vee{\mathbf{LT}}{t\in\mathcal{T}_A}[\phi[(\overline{t'},t)/(\overline{x'},x)]]
\]
in the existential case. This gives the quantifier cases by noting that 
\[
\left(\iota\tfrac{\overline{t}}{\overline{x}}\right)\tfrac{t}{x}=\iota\tfrac{(\overline{t'},t)}{(\overline{x'},x)}.
\]
\end{proof}
This immediately yields completeness theorems for $\mathcal{G}_{\omega_1}(\mathcal{L}_A)$ or respectively $\mathcal{G}_{\omega_1,\omega}(\mathcal{L}_A)$ w.r.t. L-algebras:
\begin{theorem}
Let $\mathcal{L}_A$ be any fragment of $\mathcal{L}_{\omega_1}$ or $\mathcal{L}_{\omega_1,\omega}$. For any $\Gamma\cup\{\phi\}\subseteq\mathcal{L}_A$ (with $\Gamma$ closed in the first-order case), the following are equivalent:
\begin{enumerate}
\item $\Gamma\vdash_{\mathcal{G}(\mathcal{L}_A)}\phi$;
\item $\Gamma\models_{\mathsf{L}(\mathcal{L}_A)}\phi$.
\end{enumerate}
Here, we write $\mathcal{G}$ for either $\mathcal{G}_{\omega_1}$ or $\mathcal{G}_{\omega_1,\omega}$, respectively.
\end{theorem}
\begin{proof}
``(1) implies (2)" is contained in Lemma \ref{lem:sound} or Lemma \ref{lem:fosound}, respectively. For the converse, suppose $\Gamma\not\vdash_{\mathcal{G}(\mathcal{L}_A)}\phi$ and construct the Lindenbaum-Tarski algebra $\mathbf{LT}^\Gamma$ as indicated above. 

Then Lemma \ref{lem:ltwelldef} gives that $\iota$ is a well-defined $\mathbf{LT}^\Gamma$-valued $\mathcal{L}_A$-evaluation in the propositional case and by construction, we have $\iota[\Gamma]\subseteq\{1^\mathbf{LT}\}$ but $\iota(\phi)\neq1^\mathbf{LT}$. Therefore $\Gamma\not\models_{\mathsf{L}(\mathcal{L}_A)}\phi$ as $\mathbf{LT}^\Gamma$ is an L-algebra.

In the first-order case, we construct the Lindenbaum-Tarski model $\mathfrak{LT}^\Gamma$ as above which is suitable w.r.t. $\iota$ by Lemma \ref{lem:ltsuitable}. Again by Lemma \ref{lem:ltwelldef}, the underlying algebra $\mathbf{LT}^\Gamma$ is a well-defined L-algebra and for the corresponding $\iota$, we have
\[
\overline{(\mathfrak{LT}^\Gamma,\iota)}[\Gamma]\subseteq 1^\mathbf{LT}\text{ but }\overline{(\mathfrak{LT}^\Gamma,\iota)}(\phi)=[\phi]\neq 1^\mathbf{LT}
\]
by Lemma \ref{lem:ltsuitable} and therefore $\Gamma\not\models_{\mathsf{L}(\mathcal{L}_A)}\phi$.
\end{proof}
However, restricting to countable fragments even yields a further completeness result for $\mathcal{G}^D$ w.r.t. countable chains and thus, by Lemma \ref{lem:qembed}, for $[0,1]$ which is what we will outline next.
\begin{theorem}\label{thm:distcompthmalgfrag}
Let $\mathcal{L}_A$ be a countable fragment of $\mathcal{L}_{\omega_1}$ or $\mathcal{L}_{\omega_1,\omega}$. For any $\Gamma\cup\{\phi\}\subseteq\mathcal{L}_A$ (with $\Gamma$ closed in the first-order case), the following are equivalent:
\begin{enumerate}
\item $\Gamma\vdash_{\mathcal{G}^D(\mathcal{L}_A)}\phi$;
\item $\Gamma\models_{\mathsf{DL}(\mathcal{L}_A)}\phi$;
\item $\Gamma\models_{\mathsf{C}(\mathcal{L}_A)}\phi$;
\item $\Gamma\models_{\mathsf{CC}(\mathcal{L}_A)}\phi$;
\item $\Gamma\models_{\mathbf{[0,1]_\mathbb{Q}}(\mathcal{L}_A)}\phi$;
\item $\Gamma\models_{\mathbf{[0,1]_\mathbb{R}}(\mathcal{L}_A)}\phi$.
\end{enumerate}
Here, we again write $\mathcal{G}^D$ for either $\mathcal{G}^D_{\omega_1}$ or $\mathcal{G}^D_{\omega_1,\omega}$, respectively.
\end{theorem}
\begin{proof}
``(1) implies (2)" is contained in Lemma \ref{lem:sound} or Lemma \ref{lem:fosound}, respectively . ``(2) implies (3)" follows from the fact that every chain is a distributive L-algebra. ``(3) implies (4)" and ``(4) implies (5)" as well as ``(6) implies (5)" and ``(1) implies (6)" are also immediate.\\

We thus only show ``(5) implies (1)" and for that, suppose $\Gamma\not\vdash_{\mathcal{G}^D(\mathcal{L}_A)}\phi$. Also here, we construct the corresponding Lindenbaum-Tarski algebra $\mathbf{LT}^\Gamma$ and naturally have $[\phi]\neq1^\mathbf{LT}$. Lemma \ref{lem:ltwelldef} again guarantees that $\mathbf{LT}^\Gamma$ is a well-defined L-algebra and that $\iota$ is a well-defined evaluation in the propositional case.

In the first-order case, we again construct the Lindenbaum-Tarski model $\mathfrak{LT}^\Gamma$ over the algebra $\mathbf{LT}^\Gamma$ and by Lemma \ref{lem:ltsuitable}, the model $\mathfrak{LT}^\Gamma$ is suitable w.r.t. $\iota$. With $\iota$ as the canonical variable assignment, we also get
\[
\overline{(\mathfrak{LT}^\Gamma,\iota)}(\psi) = [\psi]
\]
by the same result.

Lemma \ref{lem:ltwelldef} now yields that any of the meets
\[
\Wedge{\mathbf{LT}}{i\in\omega}([\chi]\rightarrow^\mathbf{LT}[\phi_i])=[\chi]\rightarrow^\mathbf{LT}\Wedge{\mathbf{LT}}{i\in\omega}[\phi_i]\text{ and }\Wedge{\mathbf{LT}}{i\in\omega}([\psi_i]\rightarrow^\mathbf{LT}[\chi])=\Vee{\mathbf{LT}}{i\in\omega}[\psi_i]\rightarrow^\mathbf{LT}[\chi]
\]
and additionally, in the first-order case, any of the meets
\[
\Wedge{\mathbf{LT}}{t\in\mathcal{T}_A}([\chi]\rightarrow^\mathbf{LT}[\phi[t/x]])=[\chi]\rightarrow^\mathbf{LT}[\forall x\phi]\text{ and }\Wedge{\mathbf{LT}}{t\in\mathcal{T}_A}([\phi[t/x]]\rightarrow^\mathbf{LT}[\chi])=[\exists x\phi]\rightarrow^\mathbf{LT}[\chi]
\]
are distributive for any $\chi,\bigwedge_{i\in\omega}\phi_i,\bigvee_{i\in\omega}\psi_i,\phi\in\mathcal{L}_A$ and $x\in Var_A$.

As $\mathcal{L}_A$ is countable, we can enumerate all of the above (distributive) meets and by Lemma \ref{lem:filterlemma}, there is a prime filter $F$ with $[\phi]\not\in F$ and such that $F$ preserves all the above meets (depending on the propositional or first-order case).

As $[\chi]$ reaches every element of $\mathbf{LT}^\Gamma$, Lemma \ref{lem:meetjoinfilter} gives that the map 
\[
p_F:\mathbf{LT}^\Gamma\to\mathbf{LT}^\Gamma/F, [\chi]\mapsto[\chi]_F:=[[\chi]]_F
\]
is a homomorphism of Heyting algebras with preserves all the meets/joins 
\[
\Wedge{\mathbf{LT}}{i\in\omega}[\phi_i]\text{ and }\Vee{\mathbf{LT}}{i\in\omega}[\psi_i]
\]
as well as
\[
\Wedge{\mathbf{LT}}{t\in\mathcal{T}_A}[\phi[t/x]]\text{ and }\Vee{\mathbf{LT}}{t\in\mathcal{T}_A}[\phi[t/x]]
\]
for any $\bigwedge_{i\in\omega}\phi_i,\bigvee_{i\in\omega}\psi_i,\phi\in\mathcal{L}_A$ and $x\in Var_A$. Further, $\mathbf{LT}^\Gamma/F$ is a chain as $F$ is a prime filter and $\mathbf{LT}^\Gamma$ is an L-algebra. 

In the propositional case, Lemma \ref{lem:homprop} implies that $p_F\circ\iota$ is a well-defined $\mathcal{L}_A$-evaluation into $\mathbf{LT}^\Gamma/F$ with $p_F(\iota(\phi))=[\phi]_F\neq 1^\mathbf{LT/F}$ as $[\phi]\not\in F$.

In the first-order case, since $\phi[\overline{t}/\overline{x}]\in\mathcal{L}_A$ for any $\phi\in\mathcal{L}_A$, any $\overline{t}\in(\mathcal{T}_A)^n$ and any $\overline{x}\in (Var_A)^n$ with $i_j\neq i_k$ for $j\neq k$, we in particular have that the map $p_F$ preserves the meets and joins
\[
\Wedge{\mathbf{LT}}{i\in\omega}[\phi_i[\overline{t}/\overline{x}]]\text{ and }\Vee{\mathbf{LT}}{i\in\omega}[\psi_i[\overline{t}/\overline{x}]]
\]
as well as
\[
\Wedge{\mathbf{LT}}{t\in\mathcal{T}_A}[\phi[\overline{t}/\overline{x}][t/x]]\text{ and }\Vee{\mathbf{LT}}{t\in\mathcal{T}_A}[\phi[\overline{t}/\overline{x}][t/x]].
\]
Since by Lemma \ref{lem:ltsuitable}, we have $\overline{(\mathfrak{LT}^\Gamma,\iota\tfrac{\overline{t}}{\overline{x}})}(\phi)=[\phi[\overline{t}/\overline{x}]]$, the map $p_F$ fulfills the premises of Lemma \ref{lem:fohomprop} and therefore, the model $p_F\circ\mathfrak{LT}^\Gamma$ is suitable for $\iota$ with $\overline{(p_F\circ\mathfrak{LT}^\Gamma,\iota\tfrac{\overline{t}}{\overline{x}})}(\phi)=[\phi[\overline{t}/\overline{x}]]_F$.\\

Lemma \ref{lem:qembed} guarantees the existence of an injective Heyting algebra homomorphism
\[
q:\mathbf{LT}^\Gamma/F\to\mathbf{[0,1]_\mathbb{Q}}
\]
which preserves all the meets and joins existing in $\mathbf{LT}^\Gamma/F$.

In the propositional case, by Lemma \ref{lem:homprop} we have that $q\circ (p_F\circ \iota)$ is a well-defined $\mathcal{L}_A$-evaluation with $(q\circ (p_F\circ \iota))(\phi)\neq 1$ by injectivity. Naturally, we have $(q\circ (p_F\circ \iota))[\Gamma]\subseteq\{1\}$ and therefore $\Gamma\not\models_{\mathbf{[0,1]_\mathbb{Q}}(\mathcal{L}_A)}\phi$.

In the first-order case, since $p_F\circ\mathfrak{LT}^\Gamma$ is suitable for $\iota$, this implies that in particular the conditions of Lemma \ref{lem:fohomprop} are met again and therefore, the model
\[
q\circ(p_F\circ\mathfrak{LT}^\Gamma)
\]
is suitable for $\iota$ as well. In particular, we have
\[
\overline{(q\circ(p_F\circ\mathfrak{LT}^\Gamma),\iota)}(\phi)=q([\phi]_F)\neq 1
\]
as $q$ is injective and as $[\phi]\not\in F$, i.e. $[\phi]_F\neq 1^{\mathbf{LT}^\Gamma/F}$.
\end{proof}
\begin{corollary}\label{cor:fullcompgld}
Let $\mathcal{L}_A$ be an arbitrary fragment of $\mathcal{L}_{\omega_1}$ or $\mathcal{L}_{\omega_1,\omega}$ but let $\Gamma\cup\{\phi\}\subseteq\mathcal{L}_A$ be countable with $\Gamma$ closed in the first-order case. Then, the following are equivalent:
\begin{enumerate}
\item $\Gamma\vdash_{\mathcal{G}^D(\mathcal{L}_A)}\phi$;
\item $\Gamma\models_{\mathsf{DL}(\mathcal{L}_A)}\phi$;
\item $\Gamma\models_{\mathsf{C}(\mathcal{L}_A)}\phi$;
\item $\Gamma\models_{\mathbf{[0,1]_\mathbb{R}}(\mathcal{L}_A)}\phi$.
\end{enumerate}
Here, we again write $\mathcal{G}^D$ for either $\mathcal{G}^D_{\omega_1}$ or $\mathcal{G}^D_{\omega_1,\omega}$, respectively.
\end{corollary}
\begin{proof}
The directions ``(1) implies (2)",  ``(2) implies (3)" as well as  ``(3) implies (4)" follow as before. Suppose $\Gamma\not\vdash_{\mathcal{G}^D(\mathcal{L}_A)}\phi$. Then, also
\[
\Gamma\not\vdash_{\mathcal{G}^D(\mathcal{L}_B)}\phi
\]
for $\mathcal{L}_B:=\mathrm{frag}(\Gamma\cup\{\phi\})$ since $\mathcal{L}_B\subseteq\mathcal{L}_A$. Since $\Gamma\cup\{\phi\}$ is countable, we have that $\mathcal{L}_B$ is countable as well. Theorem \ref{thm:distcompthmalgfrag} gives $\Gamma\not\models_{\mathbf{[0,1]_\mathbb{R}}(\mathcal{L}_B)}\phi$.

Therefore, in the propositional case there is an $\mathcal{L}_B$-evaluation $v:\mathcal{L}_B\to [0,1]$ such that $v[\Gamma]\subseteq\{1\}$ but $v(\phi)\neq1$. Similarly, in the first-order case, there is an interpretation $\mathfrak{I}=(\mathfrak{M},v)$ such that $\overline{\mathfrak{I}}[\Gamma]\subseteq\{1\}$ but $\overline{\mathfrak{I}}(\phi)=1$.

In $\mathbf{[0,1]_\mathbb{R}}$ (by completeness), every function $v:Var_A\to [0,1]$ has a single extension to an evaluation which we denote here by $\overline{v}$. If we set 
\[
v':p\mapsto\begin{cases}v(p)&\text{if }p\in\mathcal{L}_B,\\0&\text{otherwise},\end{cases}
\]
for $p\in Var_A$, then it is easy to see that $\overline{v}'(\psi)=v(\psi)$ for $\psi\in\mathcal{L}_B$ (as $v$ is an $\mathcal{L}_B$-evaluation) and therefore, we have $\overline{v}'[\Gamma]=1$ but $\overline{v}'(\phi)\neq 1$.

On the first-order side, any model over $\mathbf{[0,1]_\mathbb{R}}$ (again by completeness) is suitable for any variable assignment. We define $\mathfrak{M}'$ by $f^{\mathfrak{M}'}:=f^\mathfrak{M}$ and $P^{\mathfrak{M}'}:=P^\mathfrak{M}$ for function symbols $f$ and predicate symbols $P$ in $\sigma_B$ and otherwise set $f^{\mathfrak{M}'}$ and $P^{\mathfrak{M}'}$ arbitrary for symbols from $\sigma_A\setminus\sigma_B$. Further, we define $v'(x):=v(x)$ if $x\in Var_B$ and set it arbitrary otherwise for $x\in Var_A\setminus Var_B$. Again, a simple induction shows that
\[
\overline{(\mathfrak{M}',v')}(\psi)=\overline{(\mathfrak{M},v)}(\psi)
\]
for any $\psi\in\mathcal{L}_B$ which results in $\Gamma\not\models_{\mathbf{[0,1]_\mathbb{R}}(\mathcal{L}_A)}\phi$.
\end{proof}
In particular, for countable $\Gamma$ (closed, in the first-order case) in the appropriate language, we have
\[
\Gamma\vdash_{\mathcal{G}^D_{\omega_1}}\phi\text{ iff }\Gamma\models_{\mathsf{G}_{\omega_1}}\phi.
\]
as well as
\[
\Gamma\vdash_{\mathcal{G}^D_{\omega_1,\omega}}\phi\text{ iff }\Gamma\models_{\mathsf{G}_{\omega_1,\omega}}\phi.
\]
We will actually see later on that the requirement that the set of premises is countable can not be removed.
\begin{remark}\label{rem:closednotneeded}
In the first-order cases, the formulations of the various completeness results require the set of premises $\Gamma$ to be closed. This is, however, \emph{only} needed for the soundness result as can be seen by inspecting the various proofs for the converse direction: both the Lindenbaum-Tarski algebra and model do not rely on $\Gamma$ to be closed and neither do any remaining parts of the above presented completeness proofs.
\end{remark}
\section{Hypersequent Calculi for $\mathcal{G}^D_{\omega_1}$ and $\mathcal{G}^D_{\omega_1,\omega}$ and Cut-Elimination}
We now want to address structural proof theory for infinitary G\"odel logics, both for the propositional and first-order instances with operations of countable length. For that, we lift the usual approach towards structural proof theory for G\"odel logics via hypersequent calculi (see \cite{BCF2003}) to the infinitary case and provide cut-elimination theorems.
\subsection{Sequents, Hypersequents and Related Notions}
Hypersequents, as introduced by Avron \cite{Avr1996}, are multisets of sequents and the rules operating on hypersequents then allow for parallel modification of and for ``exchange of information" between sequents (as in particular exemplified by the rule $\mathsf{(com)}$ given by Avron). More formally, let $\mathcal{L}$ be either $\mathcal{L}_{\omega_1}$ or $\mathcal{L}_{\omega_1,\omega}$. A \emph{sequent} is a pair
\[
\Gamma\rhd\Delta
\]
of finite multisets $\Gamma$ and $\Delta$ of $\mathcal{L}$-formulas where $\Delta$ contains at most one element. We write $\Gamma,\Delta$ for multiset union and denote the multiset formed by (possibly equal) formulas $\phi_1,\dots,\phi_n$ by $[\phi_1,\dots,\phi_n]$. We simply write $\phi$ for $[\phi]$. A \emph{hypersequent} is then a multiset of sequents $\Gamma_i\rhd\Delta_i$ ($1\leq i\leq n$) which we denote by
\[
\Gamma_1\rhd\Delta_1\mid\dots\mid\Gamma_n\rhd\Delta_n.
\]
In any other way, we follow the notational conventions of \cite{BCF2003}. In particular, there is a \emph{canonical interpretation} $\mathcal{I}$ of hypersequents as $\mathcal{L}$-formulas: set $\mathcal{I}(\Gamma\rhd\Delta):=\bigwedge\Gamma\rightarrow\bigvee\Delta$ where $\bigwedge\Gamma$ ($\bigvee\Delta$) is the conjunction (disjunction) over all members of $\Gamma$ ($\Delta$), with the convention that $\bigwedge\emptyset:=\top$ and $\bigvee\emptyset:=\bot$. This $\mathcal{I}$ extends to hypersequents by
\[
\mathcal{I}(\Gamma_1\rhd\Delta_1\mid\dots\mid\Gamma_n\rhd\Delta_n):=\bigvee_{i=1}^n\mathcal{I}(\Gamma_i\rhd\Delta_i).
\]
The hypersequent systems which we consider are based on a hypersequent calculus for finitary propositional G\"odel logic introduced by Avron \cite{Avr1996}. This calculus naturally extends the usual sequent calculus for intuitionistic logic, lifted to the hypersequent setting, by a specific rule emulating the prelinearity axiom scheme $(\phi\rightarrow\psi)\lor (\psi\rightarrow\psi)$. Avron's calculus was extended to first-order G\"odel logics by Baaz and Zach in \cite{BZ2000}.
\subsection{The Systems $\mathcal{HG}^D_{\omega_1}$ and $\mathcal{HG}^D_{\omega_1,\omega}$}
The range of constituting rules for the various hypersequent calculi can be seen in Figure \ref{fig:hypersequentcalc}. They consist of the rules given in \cite{BC2002} for a hypersequent calculus for first-order G\"odel logic extended by four additional infinitary rules. Following Baaz and Ciabattoni \cite{BC2002}, the version of $\mathsf{(com)}$ given here differs from the one usually given (see e.g. \cite{BZ2000}) and (similar to \cite{BC2002}) this serves some technical purposes in the following cut-elimination proof. For further context, see in particular Remark 2 from \cite{BC2002}. Principal formulas are defined as commonly done.
\begin{figure}[h]
\begin{gather*}
\text{Initial Hypersequents}\\
\frac{}{\phi\rhd\phi},\phi\neq\bot\text{ and atomic}\;\mathsf{(id)}\qquad \frac{}{\bot\rhd}\;\mathsf{(\bot)}\\
\text{Structural Rules}\\
\frac{G}{G\mid\Gamma\rhd\Delta}\;\mathsf{(EW)}\qquad \frac{G\mid \Gamma\rhd\Delta\mid \Gamma\rhd\Delta}{G\mid \Gamma\rhd\Delta}\;\mathsf{(EC)}\qquad \frac{G\mid\Gamma_1,\Gamma_2\rhd\Delta_1\quad G\mid\Gamma_1,\Gamma_2\rhd\Delta_2}{G\mid\Gamma_1\rhd\Delta_1\mid\Gamma_2\rhd\Delta_2}\;\mathsf{(com)}\\
\qquad\frac{G\mid\Gamma\rhd\Delta}{G\mid\Gamma,\phi\rhd\Delta}\;\mathsf{(w,l)}\qquad \frac{G\mid\Gamma\,\rhd}{G\mid\Gamma\rhd\phi}\;\mathsf{(w,r)}\qquad \frac{G\mid\Gamma,\phi,\phi\rhd\Delta}{G\mid\Gamma,\phi\rhd\Delta}\;\mathsf{(c,l)}\\
\frac{G\mid\Gamma,\phi\rhd\Delta\quad G\mid\Gamma\rhd\phi}{G\mid\Gamma\rhd\Delta}\;\mathsf{(\mathsf{cut})}\\
\text{Logical Rules}\\
\frac{G\mid\Gamma\rhd\phi\quad G\mid\Gamma,\psi\rhd\Delta}{G\mid\Gamma,\phi\rightarrow\psi\rhd\Delta}\;\mathsf{(\rightarrow,l)}\qquad \frac{G\mid\Gamma,\phi\rhd\psi}{G\mid\Gamma\rhd\phi\rightarrow\psi}\;\mathsf{(\rightarrow,r)}\\
\frac{G\mid\Gamma,\phi_i\rhd\Delta}{G\mid\Gamma,\phi_0\land\phi_1\rhd\Delta}\;\mathsf{(\land_i,l)_{i=0,1}}\qquad \frac{G\mid\Gamma\rhd\phi\quad G\mid\Gamma\rhd\psi}{G\mid\Gamma\rhd\phi\land\psi}\;\mathsf{(\land,r)}\\
\frac{G\mid\Gamma,\phi\rhd\Delta\quad G\mid\Gamma,\psi\rhd\Delta}{G\mid\Gamma,\phi\lor\psi\rhd\Delta}\;\mathsf{(\lor,l)}\qquad
\frac{G\mid\Gamma\rhd\phi_i}{G\mid\Gamma\rhd\phi_0\lor\phi_1}\;\mathsf{(\lor_i,r)_{i=0,1}}\\
\text{Quantifier Rules}\\
\frac{G\mid\phi[t/x],\Gamma\rhd\Delta}{G\mid\forall x\phi(x),\Gamma\rhd\Delta}\;\mathsf{(\forall,l)}\qquad\frac{G\mid\Gamma\rhd\phi[a/x]}{G\mid\Gamma\rhd\forall x\phi(x)}\;\mathsf{(\forall,r)}\\\frac{G\mid\phi[a/x],\Gamma\rhd\Delta}{G\mid\exists x\phi(x),\Gamma\rhd\Delta}\;\mathsf{(\exists,l)}\qquad\frac{G\mid\Gamma\rhd\phi[t/x]}{G\mid\Gamma\rhd\exists x\phi(x)}\;\mathsf{(\exists,r)}\\
\text{Infinitary Rules}\\
\frac{G\mid\phi_j,\Gamma\rhd\Delta}{G\mid\bigwedge_{i\in\omega}\phi_i,\Gamma\rhd\Delta}\;\mathsf{(\bigwedge,l)}\qquad\frac{G\mid\Gamma\rhd\phi_i\;(i\in\omega)}{G\mid\Gamma\rhd\bigwedge_{i\in\omega}\phi_i}\;\mathsf{(\bigwedge,r)}\\\frac{G\mid\phi_i,\Gamma\rhd\Delta\;(i\in\omega)}{G\mid\bigvee_{i\in\omega} \phi_i,\Gamma\rhd\Delta}\;\mathsf{(\bigvee,l)}\qquad\frac{G\mid\Gamma\rhd\phi_j}{G\mid\Gamma\rhd\bigvee_{i\in\omega}\phi_i}\;\mathsf{(\bigvee,r)}
\end{gather*}
\caption{The various hypersequent rules.}
\label{fig:hypersequentcalc}
\end{figure}
The rules $\mathsf{(\forall,l)}$ and $\mathsf{(\exists,r)}$ are supposed to fulfill the eigenvariable condition: the variable $a$ is free and does not occur in the lower hypersequent. We refer with $\mathcal{HG}^D_{\omega_1}$ to all initial, structural, logical and infinitary rules (over the propositional language) and with $\mathcal{HG}^D_{\omega_1,\omega}$ to $\mathcal{HG}^D_{\omega_1}$ (now over the first-order language) extended with the quantifier rules.\\

Deductions in the hypersequent calculi are defined by natural infinitary generalizations of the usual definition: deductions are countable (possibly infinite) rooted trees where every node is labeled with a hypersequent and every edge is labeled with a rule such that the arities of the rules are respected and the applications are correct. If $d$ is such derivation with a root hypersequent $H$, then we write $d\vdash_{\mathcal{HG}^D_{\omega_1,\omega}}H$ or $d\vdash_{\mathcal{HG}^D_{\omega_1}}H$ (depending on the used language and systems). We omit the system if it is arbitrary or clear from the context and then write $d\vdash H$. We write $\vdash H$ if there is any derivation $d$ with $d\vdash H$.\\

One quickly verifies that the non-atomic version 
\[
\frac{}{\phi\rhd\phi}
\]
of $\mathsf{(id)}$ for \emph{arbitrary} $\phi$ is derivable in the systems and that they are complete w.r.t. to the Hilbert-type systems introduced before in the following sense:
\begin{theorem}
For any $\phi\in\mathcal{L}_{\omega_1}$, $\vdash_{\mathcal{HG}^D_{\omega_1}}\rhd\phi$ if, and only if $\vdash_{\mathcal{G}^D_{\omega_1}}\phi$. Similarly for $\mathcal{L}_{\omega_1,\omega}$, $\mathcal{HG}^D_{\omega_1,\omega}$ and $\mathcal{G}^D_{\omega_1,\omega}$.
\end{theorem}
\begin{proof}
The only thing we want to remark is that the infinitary distributivity axiom
\[
\tag{D}\bigwedge_{i\in\omega}(\phi\lor\psi_i)\rightarrow\left(\phi\lor\bigwedge_{i\in\omega}\psi_i\right).
\]
is derivable in the calculus, akin to the analogous derivation of the first order ($QS$) axiom e.g. given in \cite{BCF2003}. Concretely, we have the following derivations
\begin{prooftree}
\AxiomC{$\phi\rhd\phi$}
\AxiomC{$\psi_i\rhd\psi_i$}
\RightLabel{$(\mathsf{com})$}
\BinaryInfC{$\phi\rhd\psi_i\mid\psi_i\rhd\phi$}
\AxiomC{$\phi\rhd\phi$}
\RightLabel{$(\mathsf{\lor,l})$}
\BinaryInfC{$\phi\rhd\psi_i\mid\phi\lor\psi_i\rhd\phi$}
\AxiomC{$\psi_i\rhd\psi_i$}
\RightLabel{$(\mathsf{\lor,l})$}
\BinaryInfC{$\phi\lor\psi_i\rhd\psi_i\mid\phi\lor\psi_i\rhd\phi$}
\RightLabel{$2x(\mathsf{\bigwedge,l})$}
\UnaryInfC{$\bigwedge_{i\in\omega}(\phi\lor\psi_i)\rhd\psi_i\mid\bigwedge_{i\in\omega}(\phi\lor\psi_i)\rhd\phi$}
\end{prooftree}
for any $i\in\omega$. Using those as input for $\mathsf{(\bigwedge,r)}$, we get
\begin{prooftree}
\AxiomC{$\vdots$}
\RightLabel{$(\mathsf{\bigwedge,r})$}
\UnaryInfC{$\bigwedge_{i\in\omega}(\phi\lor\psi_i)\rhd\bigwedge_{i\in\omega}\psi_i\mid\bigwedge_{i\in\omega}(\phi\lor\psi_i)\rhd\phi$}
\RightLabel{$(\mathsf{\lor_0,r}),(\mathsf{\lor_1,r})$}
\UnaryInfC{$\bigwedge_{i\in\omega}(\phi\lor\psi_i)\rhd\phi\lor\bigwedge_{i\in\omega}\psi_i\mid\bigwedge_{i\in\omega}(\phi\lor\psi_i)\rhd\phi\lor\bigwedge_{i\in\omega}\psi_i$}
\RightLabel{$(\mathsf{EC})$}
\UnaryInfC{$\bigwedge_{i\in\omega}(\phi\lor\psi_i)\rhd\phi\lor\bigwedge_{i\in\omega}\psi_i$}
\RightLabel{$(\mathsf{\rightarrow,r}).$}
\UnaryInfC{$\rhd\bigwedge_{i\in\omega}(\phi\lor\psi_i)\rightarrow\phi\lor\bigwedge_{i\in\omega}\psi_i$}
\end{prooftree}
\end{proof}
For the upcoming proof of the cut-elimination theorem, we introduce versions of the calculi given above which use \emph{sets of formulas} for sequents and \emph{sets of sequents} for hypersequents. With this set-version, we follow both Tait \cite{Tai1968} as well as Baaz and Ciabattoni \cite{BC2002} and the advantage is that we can omit both the internal and external contraction rules in the set-version. We still write  $\Gamma\rhd\Delta$ for set-sequents where $\Delta$ is still at most a singleton and in this context write $\Gamma,\Delta$ for the ordinary union of $\Gamma$ and $\Delta$. To emphasize the set-version of the hypersequents, we denote the respective objects by $G\cup\{\Gamma\rhd\Delta\}$. 

All the previously introduced rules (besides external and internal contraction) can be naturally reformulated using the set-hypersequents and we use the same naming in theses cases. We write $\mathcal{SHG}^D_{\omega_1}$ and $\mathcal{SHG}^D_{\omega_1,\omega}$ for these set-versions of the previous calculi. We use the same notion of proof defined via countable trees and write $d\vdash^s_{\mathcal{SHG}^D_{\omega_1,\omega}} H$ or $\vdash^s_{\mathcal{SHG}^D_{\omega_1}}H$ if $d$ is a derivation tree with $H$ as a root set-hypersequent in the respective systems. Again, if we don't want to be specific about the system, we write $d\vdash^s H$. Also, derivations $d\vdash^s H$ with additional set-hypersequent assumptions are defined as always. The notion of substitution naturally carries over to sequents, (set-)hypersequents and derivations and we use the same notation as in formulas and terms.
\subsection{Cut-Elimination in the Sch\"utte-Tait Style}
We now turn to cut-elimination. Most cut-elimination methods fall into one of two categories: Gentzen style methods \cite{Gen1934} which remove highest cuts and Sch\"utte-Tait style methods \cite{Sch1960,Tai1968} which remove most complex cuts (in the sense of occurring logical symbols or a similar complexity measure). The notion of highest cut does not generally result in terminating procedures with systems which have infinitary rules and it is thus not surprising that we opt for a Sch\"utte-Tait style proof.

For proving cut-elimination, we closely follow the argument given in \cite{BC2002} by Baaz and Ciabattoni where the authors provide a Sch\"utte-Tait style cut-elimination proof for the calculus which we have used (together with its propositional fragment) as a basis for our infinitary extensions. Naturally, the finitary notions used there have to be appropriately extended to arbitrary countable ordinals and we do this in a similar vein as in Tait's work \cite{Tai1968}.\\

Before moving on to the technical results, we will need to introduce various measures on proofs and formulas. For this, we first give a short primer on the relevant notions regarding ordinals. For a general overview over ordinal arithmetic and further notions, see \cite{Sie1958}. We write $+$ and $\cdot$ for the usual ordinal addition and multiplication. Given a family of ordinals $\alpha_i$, we write $\sup^+_i\alpha_i$ for the smallest ordinal greater than every $\alpha_i$ and $\sup_i\alpha_i$ for the smallest ordinal greater or equal than every $\alpha_i$.

Further, we will need the \emph{natural sum} and \emph{multiplication} on ordinals, also called the Hessenberg sum and multiplication (see \cite{Sie1958}), which we denote by $x\oplus y$ and $x\otimes y$, respectively. The precise definition for $\oplus$ and $\otimes$ can be given using Cantor normal forms, among others, but we will only need certain properties of $\oplus$, $\otimes$ and their interplay with $\sup^+$ which we list in the following remark.
\begin{remark}
$\oplus$ and $\otimes$ are commutative, associative and monotone in both arguments. Further, we have $\alpha\oplus 1=\alpha+1$ for any ordinal $\alpha$ and $\alpha\otimes (\beta\oplus\gamma)$ = $(\alpha\otimes\beta)\oplus(\alpha\otimes\gamma)$. We have $\sup^{(+)}_i(\beta\oplus\alpha_i)\leq\beta\oplus\sup^{(+)}_i\alpha_i$ and $\sup^{(+)}_i(\beta\otimes\alpha_i)$ $\leq$ $\beta\otimes\sup^{(+)}_i\alpha_i$ where $\sup^{(+)}$ is either $\sup$ or $\sup^+$.
\end{remark}
There is a notion of exponentiation derived from natural multiplication, first considered by de Jongh and Parikh \cite{deJP1977}, which we denote by $\alpha^{\otimes\beta}$ and call \emph{(super-Jacobsthal) exponentiation}, following Altman \cite{Alt2017}. This exponentiation can be formally defined using transfinite recursion by 
\begin{enumerate}
\item $\alpha^{\otimes 0}:=1$ for any $\alpha$,
\item $\alpha^{\otimes (\beta+1)}:=\alpha^{\otimes\beta}\otimes\alpha$ for any $\alpha,\beta$,
\item $\alpha^{\otimes\beta}:=\sup_{\gamma <\beta}\alpha^{\otimes\gamma}$ for any $\alpha$ and any limit ordinal $\beta$,
\end{enumerate}
and has, in particular, the following properties: $\alpha^{\otimes\beta}$ is strictly increasing and continuous in $\beta$ and
\[
\alpha^{\otimes(\beta\oplus\gamma)}=\alpha^{\otimes\beta}\otimes\alpha^{\otimes\gamma}.
\]
The first two are immediate by the monotonicity of $\otimes$ and the definition of exponentiation. A proof of the latter can be found in \cite{Alt2017}.\\

Inspired by Tait \cite{Tai1968}, we define the function $\chi^0(\alpha):=4^{\otimes \alpha}$ and derived from that, we define $\chi^z$ as the function enumerating the common fixed points of $\chi^w$ for all $w<z$. These $\chi^z$ are the Veblen iterations of $\chi^0$ which exist (for all countable ordinals) by \cite{Veb1908} as $\chi^0$ is continuous and strictly increasing. Further, all $\chi^z$ are also continuous and strictly increasing.\\

We define the \emph{complexity} $\vert\phi\vert$ of $\phi\in\mathcal{L}_{\omega_1}$ recursively by
\begin{enumerate}
\item $\vert\phi\vert := 0$ for atomic $\phi$,
\item $\vert\phi\circ\psi\vert := \max\{\vert\phi\vert,\vert\psi\vert\}+1$,
\item $\vert\bigwedge_{i\in\omega}\phi_i\vert:=\vert\bigvee_{i\in\omega}\phi_i\vert:=\sup^+_{i\in\omega}\vert\phi_i\vert$.
\end{enumerate}
$\vert\cdot\vert$ is extended to $\mathcal{L}_{\omega_1,\omega}$ by adding the clause
\begin{enumerate}
\setcounter{enumi}{3}
\item $\vert Qx\phi\vert:=\vert\phi\vert+1$,
\end{enumerate}
for $Q\in\{\forall,\exists\}$.\\

Let $d$ be a derivation (of normal or set-hypersequents) with $d_i$, $i<k\leq\omega$, as its direct predecessors (i.e., those subderivations proving the assumptions of the last rule in $d$ as the \emph{direct subderivations} from \cite{Tai1968}). As natural generalizations of the notions defined in \cite{BC2002} (and in some way akin to \cite{Tai1968}) we set $\vert d\vert:=\sup^+_{i<k}\vert d_i\vert$ if the last rule was not a weakening and $\vert d\vert :=\vert d_0\vert$ otherwise. We define $w(d):=\sup_{i<k}w(d_i)$ if the last rules was not an internal weakening and $w(d):=w(d_0)+1$ otherwise. Similarly, we define $W(d)$ by using external instead of internal weakenings. 

Lastly, we recursively define $\rho(d)$ by
\begin{enumerate}
\item $\rho(d):=0$ if $d$ is cut-free,
\item $\rho(d):=\sup_{i<k}\rho(d_i)$ if the last inference is not a cut,
\item $\rho(d):=\max\{\vert\phi\vert+1,\rho(d_0),\rho(d_1)\}$ if the last inference is a cut with cut-formula $\phi$.
\end{enumerate}
$\vert d\vert$ is called the \emph{rank of $d$} and $\rho(d)$ is called the \emph{cut-degree of $d$}. All of these naturally extend to derivations with assumptions.
\begin{lemma}\label{lem:sub}
If $d\vdash^s H$, then $d[t/x]\vdash^s H[t/x]$ with $\vert d[t/x]\vert=\vert d\vert$ and $\rho(d[t/x])=\rho(d)$.
\end{lemma}
It is straightforward to check that $d[t/x]$ is a correct proof since $d$ is. Both $d[t/x]\vdash^s H[t/x]$ and the other properties are then immediate.\\

It is clear that $\vdash$ and $\vdash^s$ prove essentially the same theorems (modulo applications of contractions) and using the above notions, we can state the following result on the impact of change between $\mathcal{HG}$ and $\mathcal{SHG}$ on the rank of derivations. For that, we write $H^s$ for the set-hypersequent obtained from a hypersequent $H$ by removing all multiplicities of formulas and sequents (and treating the resulting objects as sets).

Following \cite{BC2002}, we call a hypersequent $H$ \emph{(1-1-)reduced} if no formula occurs more than once in any multiset and if no sequent occurs more than once in the hypersequent. These hypersequents can be naturally seen as set-hypersequents.
\begin{lemma}
Let $H$ be a reduced hypersequent. If $d'\vdash^s H^s$, then there is a $d\vdash H$ with $\vert d\vert\leq (2\otimes\vert d'\vert)\oplus w(d')$ and $\rho(d')=\rho(d)$. Conversely, if $d\vdash H$, then there is a $d'\vdash^s H^s$ with $\vert d'\vert\leq\vert d\vert$ and $\rho(d')=\rho(d)$.
\end{lemma}
We omit the proof as it is a natural generalization of the respective finitary result of Baaz and Ciabattoni \cite{BC2002}.\\

The proofs which we give in the following rely on a (formal) tracking of the cut formula through the proof based on so-called \emph{decorations} as introduced by Baaz and Ciabattoni in \cite{BC2002}, extended to the infinitary case. 
\begin{definition}
Given $d\vdash^s H$ and given a \emph{decoration of $H$}, that is $H$ where some (but not necessarily all) occurrences of a formula $\phi$ are decorated, denoted by $\phi^*$, the decorated version of $d$ is defined by recursion on the tree: if we have a decoration of an occurring hypersequent $H'$, then the premises are decorated according to which rule was used to derive $H'$. The definitions for the rules $\mathsf{(EW)}$, $\mathsf{(w,l)}$, $\mathsf{(w,r)}$, $\mathsf{(com)}$ and $\mathsf{(com)}$ are exactly as in \cite{BC2002}. Suppose the rule used is a (possibly infinitary) logical rule with arity $k\leq\omega$, i.e. we have
\[
\frac{G\cup\{\Gamma_i\rhd\Delta_i\mid i<k\}}{G\cup\{\Gamma\rhd\Delta\}}.
\]
Then if
\begin{enumerate}
\item $\phi$ is the principal formula of the rule: if $\phi^*\in\Gamma$ then $\phi$ is decorated in $\Gamma_i$ iff $\phi$ already occurs in $\Gamma_i$.
\item $\phi$ is not the principal formula of the rule then $\phi$ is decorated in $\Gamma_i$ or $\Delta_i$ iff it is decorated in $\Gamma$ or $\Delta$, respectively. 
\end{enumerate}
Further, in both cases, $G\setminus\{\Gamma_i\rhd\Delta_i\mid i<k\}$ is decorated as in the conclusion.
\end{definition}
\begin{lemma}\label{lem:inversion}
The following inversions are valid:
\begin{enumerate}
\item[(i)] If $d\vdash^s G\cup\{\Gamma,\phi\lor\psi\rhd\Delta\}$, then there are proofs $d_0\vdash^s G\cup\{\Gamma,\phi\rhd\Delta\}$ and $d_1\vdash^s G\cup\{\Gamma,\psi\rhd\Delta\}$.
\item[(ii)] If $d\vdash^s G\cup\{\Gamma,\phi\land\psi\rhd\Delta\}$, then there is a proof $d_0\vdash^s G\cup\{\Gamma,\phi,\psi\rhd\Delta\}$.
\item[(iii)] If $d\vdash^s G\cup\{\Gamma\rhd\phi\land\psi\}$, then there are proofs $d_0\vdash^s G\cup\{\Gamma\rhd\phi\}$ and $d_1\vdash^s G\cup\{\Gamma\rhd\psi\}$.
\item[(iv)] If $d\vdash^s G\cup\{\Gamma\rhd\phi\rightarrow\psi\}$, then there is a proof $d_0\vdash^s G\cup\{\Gamma,\phi\rhd\psi\}$.
\item[(v)] If $d\vdash^s G\cup\{\Gamma,\exists x\phi\rhd\Delta\}$, then there is a proof $d_0\vdash^s G\cup\{\Gamma,\phi[a/x]\rhd\Delta\}$.
\item[(vi)] If $d\vdash^s G\cup\{\Gamma\rhd\forall x\phi\}$, then there is a proof $d_0\vdash^s G\cup\{\Gamma\rhd\phi[a/x]\}$.
\item[(vii)] If $d\vdash^s G\cup\{\Gamma,\bigvee_{i\in\omega}\phi_i\rhd\Delta\}$, then there are proofs $d_j\vdash^s G\cup\{\Gamma,\phi_j\rhd\Delta\}$ for each $j\in\omega$.
\item[(viii)] If $d\vdash^s G\cup\{\Gamma\rhd\bigwedge_{i\in\omega}\phi_i\}$, then there are proofs $d_j\vdash^s G\cup\{\Gamma\rhd\phi_j\}$ for each $j\in\omega$.
\end{enumerate}
In any case, we respectively have $\rho(d_i)\leq\rho (d)$ and $\vert d_i\vert\leq \vert d\vert$.
\end{lemma}
\begin{proof}
The proof for items (i) to (vi) follows exactly the reasoning of \cite{BC2002} for the respective finitary result (see Lemma 4 there). We give the proofs for (vii) and (viii) in the same spirit.
\begin{enumerate}
\item[(vii)] We consider a decoration of $d$ starting with $G\cup\{\Gamma,(\bigvee_{i\in\omega}\phi_i)^*\rhd\Delta\}$. Replace every occurring $\Gamma',(\bigvee_{i\in\omega}\phi_i)^*\rhd\Delta'$ by $\Gamma',\phi_j\rhd\Delta'$. Delete all the subderivations but the {$j$-th} one above any application of $\mathsf{(\bigvee,l)}$ where $(\bigvee_{i\in\omega}\phi_i)^*$ occurs as a decorated formula and is principle. As all initial hypersequents are atomic, correctness of the resulting $d_j$ can be shown by an induction over $\vert d\vert\oplus w(d)\oplus W(d)$. Clearly $\vert d_j\vert\leq\vert d\vert$ and $\rho(d_j)\leq \rho(d)$.
\item[(viii)] We consider a decoration of $d$ starting from $G\cup\{\Gamma\rhd(\bigwedge_{i\in\omega}\phi_i)^*\}$. Replace every occurring $\Gamma'\rhd(\bigwedge_{i\in\omega}\phi_i)^*$ by $\Gamma'\rhd\phi_j$ and delete all subderivations but the {$j$-th} one above any application of $\mathsf{(\bigwedge,r)}$ in which $(\bigwedge_{i\in\omega}\phi_i)^*$ occurs decorated and is principle. Again, the correctness of $d_j$ follows by a straightforward induction on $\vert d\vert\oplus w(d)\oplus W(d)$ and we have $\vert d_j\vert\leq \vert d\vert$ and $\rho(d_j)\leq\rho (d)$ already by construction.
\end{enumerate}
\end{proof}
\begin{lemma}\label{lem:atomicCutElim}
Suppose $d\vdash^s G\cup\{\Gamma,\phi\rhd\Delta\}$ where $\phi$ is atomic and not the cut-formula of any cut in $d$. Then for any $\Sigma$, there is a $d'$ with assumption $G\cup\{\Sigma\rhd\phi\}$ such that $d'\vdash G\cup\{\Gamma,\Sigma\rhd\Delta\}$ with $\vert d'\vert\leq\vert d\vert$ and $\rho(d')\leq\rho(d)$.
\end{lemma}
\begin{proof}
The proof can be easily obtained by generalizing the proof of the respective finitary result from \cite{BC2002}: Decorate $d$, starting from $G\cup\{\Gamma,\phi^*\rhd\Delta\}$, replace any occurring $\{\Gamma',\phi^*\rhd\Delta'\}$ by $\{\Gamma',\Sigma\rhd\Delta'\}$ and add $G$ to every set-hypersequent. This tree now needs to be corrected to yield a correct proof. As in \cite{BC2002}, weakenings which produce decorated $\phi^*$ are replaced by (potentially more) weakenings producing $\Sigma$. Initial sequents
\[
\frac{}{\psi\rhd\psi}\mathsf{(id)}\text{ and }\frac{}{\bot\rhd}\mathsf{(\bot)}
\]
which don't introduce $\phi^*$ are respectively replaced by 
\[
\frac{\frac{}{\psi\rhd\psi}}{G\cup\{\psi\rhd\psi\}}\text{ and }\frac{\frac{}{\bot\rhd}}{G\cup\{\bot\rhd\}}
\]
where we added sufficiently many weakenings to introduce $G$,
while initial sequents introducing $\phi^*$ are replaced by $G\cup\{\Sigma\rhd\phi\}$ which is an allowed assumption for $d'$. That proof is now correct as $\phi$ is not the cut-formula of any cut in $d$. As weakenings do not lengthen $\vert d\vert$, we get $\vert d'\vert\leq\vert d\vert$ and as we didn't introduce any new cut, we get $\rho(d')\leq\rho(d)$.
\end{proof}
\begin{lemma}\label{lem:reduction}
Suppose $d_0\vdash^s G\cup\{\Gamma\rhd\phi\}$ and $d_1\vdash^s G\cup\{\Gamma,\phi\rhd\Delta\}$ with $\rho(d_i)\leq\vert\phi\vert$ for $i=0,1$. Then there is a $d\vdash^s G\cup\{\Gamma\rhd\Delta\}$ with $\rho(d)\leq\vert\phi\vert$ and $\vert d\vert\leq 2\otimes(\vert d_0\vert\oplus \vert d_1\vert)$.
\end{lemma}
\begin{proof}
The proof is a natural extension of the corresponding finitary result from \cite{BC2002}. Note that since $\rho(d_i)\leq\vert\phi\vert$, $\phi$ is not the cut-formula of any cut in $d_i$.\\

For $\phi=\bot$, decorate $d_0$ starting from $G\cup\{\Gamma\rhd\bot^*\}$ and replace any $\{\Gamma'\rhd\bot^*\}$ by $\{\Gamma'\rhd\Delta\}$. We now have to correct the proof at the points where $\bot^*$ originates. In this simple case, $\bot^*$ arises by either an internal or external weakening. The internal weakenings can be either removed if $\Delta$ is empty or replaced by an internal weakening with $\Delta$ if $\Delta$ is nonempty. Similarly, the external weakenings get appropriately replaced.\\

For atomic $\phi\neq\bot$, note first that $\rho(d_i)\leq\vert\phi\vert=0$ implies that $d_1$ and $d_0$ are cut-free. In particular, $\phi$ is not the cut-formula of any cut in $d_1$ and thus, there is a derivation $d_1'$ with assumption $G\cup\{\Gamma\rhd\phi\}$ such that
\[
d_1'\vdash^s G\cup\{\Gamma\rhd\Delta\}
\]
as well as $\vert d'_1\vert\leq\vert d_1\vert$ and $\rho(d_1')=0$. We form $d$ by replacing every assumption $G\cup\{\Gamma\rhd\phi\}$ with the proof $d_0$. It is straightforward to check that $\vert d\vert\leq\vert d_0\vert\oplus\vert d_1\vert$ $\leq$ $2\otimes (\vert d_0\vert\oplus\vert d_1\vert)$ and we have $\rho(d)=0$ by Lemma \ref{lem:atomicCutElim}.\\

For $\phi=\bigvee_{i\in\omega}\phi_i$, consider a decoration of $d_0$ starting from $G\cup\{\Gamma\rhd(\bigvee_{i\in\omega}\phi_i)^*\}$. Replace any occurring $\{\Sigma\rhd(\bigvee_{i\in\omega}\phi_i)^*\}$ by $\{\Gamma,\Sigma\rhd\Delta\}$ and add $G$ to every set-hypersequent and $\Gamma$ to every premise of any sequent. This resulting tree is not a correct proof anymore and we consider the following correction steps on the initial rules and on the rules which introduce a decorated instance of $\phi$:
\begin{enumerate}
\item Replace any initial rule 
\[
\frac{}{\psi\rhd\psi}\mathsf{(id)}\text{ or }\frac{}{\bot\rhd}\mathsf{(\bot)}
\]
by
\begin{prooftree}
\AxiomC{}
\RightLabel{$\mathsf{(id)}$}
\UnaryInfC{$\psi\rhd\psi$}
\UnaryInfC{$G\cup\{\Gamma,\psi\rhd\psi\}$}
\end{prooftree}
or
\begin{prooftree}
\AxiomC{}
\RightLabel{$\mathsf{(\bot)}$}
\UnaryInfC{$\bot\rhd$}
\UnaryInfC{$G\cup\{\Gamma,\bot\rhd\}$}
\end{prooftree}
respectively, using sufficiently many applications of $\mathsf{(EW)}$.
\item Suppose $(\bigvee_{i\in\omega}\phi_i)^*$ originates as the principal formula of a logical rule. Then, replace any part of the form
\begin{prooftree}
\AxiomC{$\vdots$}
\UnaryInfC{$G'\cup\{\Gamma'\rhd\phi_j\}$}
\RightLabel{$(\mathsf{\bigvee,r})$}
\UnaryInfC{$G'\cup\{\Gamma'\rhd (\bigvee_{i\in\omega}\phi_i)^*\}$}
\end{prooftree}
occurring in the decorated version of $d_0$ by
\begin{prooftree}
\AxiomC{$\vdots$}
\UnaryInfC{$G\cup (\overline{G}'\cup\{\Gamma,\Gamma'\rhd\phi_j\})[\Delta/(\bigvee_{i\in\omega}\phi_i)^*]$}
\AxiomC{$\vdots\; d_1^{-1}$}
\UnaryInfC{$G\cup \overline{G}'\cup\{\Gamma,\Gamma',\phi_j\rhd\Delta\}$}
\RightLabel{$\mathsf{(cut)}$}
\BinaryInfC{$G\cup \overline{G}'\cup\{\Gamma,\Gamma'\rhd\Delta\}$}
\end{prooftree}
where $d_1^{-1}$ is obtained by Lemma \ref{lem:inversion} from $d_1\vdash G\cup\{\Gamma,\bigvee_{i\in\omega}\phi_i\rhd\Delta\}$. Here, $\overline{G}'$ is the set-hypersequent $G'$ after the possible internal replacements. In the above presentation, we have suppressed the needed external and internal weakenings as they have no effect on the resulting rank.
\item If $(\bigvee_{i\in\omega}\phi_i)^*$ originates from an internal weakening, replace any such
\begin{prooftree}
\AxiomC{$\vdots$}
\UnaryInfC{$G'\cup\{\Gamma'\rhd\}$}
\RightLabel{$\mathsf{(w,r)}$}
\UnaryInfC{$G'\cup\{\Gamma'\rhd(\bigvee_{i\in\omega}\phi_i)^*\}$}
\end{prooftree}
by 
\begin{prooftree}
\AxiomC{$\vdots$}
\UnaryInfC{$\overline{G}'\cup G\cup\{\Gamma'\rhd\}$}
\RightLabel{$\mathsf{(w,r)}$}
\UnaryInfC{$\overline{G}'\cup G\cup\{\Gamma,\Gamma'\rhd\Delta\}$}
\end{prooftree}
if $\Delta$ is nonempty and remove them otherwise.
\item If $(\bigvee_{i\in\omega}\phi_i)^*$ originates from an external weakening, replace any 
\begin{prooftree}
\AxiomC{$\vdots$}
\UnaryInfC{$G'$}
\RightLabel{$\mathsf{(EW)}$}
\UnaryInfC{$G'\cup\{\Gamma'\rhd(\bigvee_{i\in\omega}\phi_i)^*\}$}
\end{prooftree}
by 
\begin{prooftree}
\AxiomC{$\vdots$}
\UnaryInfC{$\overline{G}'\cup G$}
\RightLabel{$\mathsf{(EW)}.$}
\UnaryInfC{$\overline{G}'\cup G\cup\{\Gamma,\Gamma'\rhd\Delta\}$}
\end{prooftree}
\end{enumerate}
The resulting proof $d$ is correct as can be verified by transfinite induction on $\vert T\vert\oplus w(T)\oplus W(T)$ for subderivations $T$ of $d_0$. Here, it is important that $\phi$ is not the cut-formula of any cut.

Furthermore, any newly introduced cut uses one of the $\phi_i$ as a cut-formula. We thus get, using $\rho(d_1^{-1})\leq\rho(d_1)$ from Lemma \ref{lem:inversion}, that
\[
\rho(d)\leq\max\left\{\rho(d_0),\rho(d_1),\sup_{i\in\omega}(\vert\phi_i\vert+1)\right\}=\max\left\{\rho(d_0),\rho(d_1),\sideset{}{^{+}}\sup_{i\in\omega} \vert\phi_i\vert\right\}\leq\vert\phi\vert
\]
as $\rho(d_i)\leq\vert\phi\vert$ for $i=0,1$ by assumption and $\vert\phi\vert=\sup^+_{i\in\omega}\vert\phi_i\vert$ by definition.\\

Regarding the rank, let $T'$ be the replacement of any rooted sub derivation $T$ of $d_0$ \emph{after} the replacement and correction procedure. Then we can prove, by transfinite induction on $\vert T\vert\oplus w(T)\oplus W(T)$, that $\vert T'\vert\leq 2\otimes (\vert d_1^{-1}\vert\oplus\vert T\vert)$. 
\begin{enumerate}
\item If $T$ is just a single initial rule, then $\vert T'\vert=\vert T\vert$ as only weakenings were added.
\item If $T$ ends with an application of $\mathsf{(\bigvee,r)}$ which introduces $\phi^*$ in the annotation, let $T_0$ be its preceding derivation. Note that $T'$ ends with an application of cut with preceding derivations $T_0'$ and $d_1^{-1}$ (modulo additional weakenings). We get
\begin{align*}
\vert T'\vert&=\max\{\vert T_0'\vert,\vert d_1^{-1}\vert\}+1\\
&\leq \max\{2\otimes(\vert T_0\vert\oplus\vert d_1^{-1}\vert),\vert d_1^{-1}\vert\}+1\\
&\leq 2\otimes((\vert T_0\vert+1)\oplus\vert d_1^{-1}\vert)\\
&=2\otimes(\vert T\vert\oplus\vert d_1^{-1}\vert)
\end{align*}
which completes this case. Here, the first inequality follows from the induction hypothesis.
\item If $T$ ends with a weakening which introduces $\phi^*$, then we still have $\vert T_0\vert\oplus w(T_0)\oplus W(T_0)<\vert T\vert\oplus w(T)\oplus W(T)$ and thus, we get
\[
\vert T'\vert=\vert T_0'\vert\leq 2\otimes (\vert d_1^{-1}\vert\oplus\vert T_0\vert)=2\otimes (\vert d_1^{-1}\vert\oplus\vert T\vert)
\]
as by definition $\vert T_0\vert =\vert T\vert$ and as only weakenings were added.
\item If the last rule of $T$ does not introduce $\phi^*$ and is not a weakening, then for the preceding derivations $T'_i$ of $T'$, we have
\[
\vert T'\vert=\sideset{}{^{+}}\sup_{i\in\omega}\vert T'_i\vert\leq\sideset{}{^{+}}\sup_{i\in\omega}2\otimes (\vert d_1^{-1}\vert\oplus\vert T_i\vert)\leq 2\otimes (\vert d_1^{-1}\vert\oplus\sideset{}{^{+}}\sup_{i\in\omega}\vert T_i\vert)=2\otimes (\vert d_1^{-1}\vert\oplus\vert T\vert).
\]
\item If the last rules does not introduce $\phi^*$ but is a weakening, then the reasoning is as in case (3).
\end{enumerate}
We in particular thus have (taking $d=d_0'$)
\[
\vert d\vert\leq 2\otimes(\vert d_1^{-1}\vert\oplus\vert d_0\vert)\leq2\otimes(\vert d_1\vert\oplus\vert d_0\vert)
\]
as $\vert d_1\vert\leq\vert d_1^{-1}\vert$ by Lemma \ref{lem:inversion}.\\

For $\phi=\bigwedge_{i\in\omega}\phi_i$, consider a decoration of $d_1$ starting from $G\cup\{\Gamma,(\bigwedge_{i\in\omega}\phi_i)^*\rhd\Delta\}$. Replace any occurring $\{\Gamma',(\bigwedge_{i\in\omega}\phi_i)^*\rhd\Delta'\}$ by $\{\Gamma',\Gamma\rhd\Delta'\}$, add $G$ to every set-hypersequent and add $\Gamma$ to the premise of every sequent. This resulting tree is not a correct proof anymore and we consider the following correction steps:
\begin{enumerate}
\item For initial rules, this is the same as the previous correction step (1).
\item If $(\bigwedge_{i\in\omega}\phi_i)^*$ originates as the principal formula of a logical rule, then replace any part of the form
\begin{prooftree}
\AxiomC{$\vdots$}
\UnaryInfC{$G'\cup\{\Gamma',\phi_j\rhd\Delta'\}$}
\RightLabel{$(\mathsf{\bigvee,r})$}
\UnaryInfC{$G'\cup\{\Gamma',(\bigwedge_{i\in\omega}\phi_i)^*\rhd \Delta'\}$}
\end{prooftree}
of the current proof by
\begin{prooftree}
\AxiomC{$\vdots$}
\UnaryInfC{$G\cup (\overline{G}'\cup\{\Gamma,\Gamma',\phi_j\rhd\Delta'\}[\Gamma/(\bigwedge_{i\in\omega}\phi_i)^*]$}
\AxiomC{$\vdots\; d_0^{-1}$}
\UnaryInfC{$G\cup \overline{G}'\cup\{\Gamma,\Gamma'\rhd\phi_j\}$}
\RightLabel{$\mathsf{(cut)}$}
\BinaryInfC{$G\cup \overline{G}'\cup\{\Gamma,\Gamma'\rhd\Delta\}$}
\end{prooftree}
where $d_0^{-1}$ is obtained via Lemma \ref{lem:inversion} from $d_0\vdash G\cup\{\Gamma\rhd\bigwedge_{i\in\omega}\phi_i\}$.
\item If $(\bigwedge_{i\in\omega}\phi_i)^*$ originates from an internal weakening, replace any such
\begin{prooftree}
\AxiomC{$\vdots$}
\UnaryInfC{$G'\cup\{\Gamma'\rhd\Delta'\}$}
\RightLabel{$\mathsf{(w,r)}$}
\UnaryInfC{$G'\cup\{\Gamma',(\bigwedge_{i\in\omega}\phi_i)^*\rhd\Delta'\}$}
\end{prooftree}
by 
\begin{prooftree}
\AxiomC{$\vdots$}
\UnaryInfC{$\overline{G}'\cup G\cup\{\Gamma'\rhd\Delta'\}$}
\UnaryInfC{$\overline{G}'\cup G\cup\{\Gamma',\Gamma\rhd\Delta'\}$}
\end{prooftree}
using stepwise internal weakenings for members of $\Gamma$.
\item If $(\bigwedge_{i\in\omega}\phi_i)^*$ originates from an external weakening, replace any 
\begin{prooftree}
\AxiomC{$\vdots$}
\UnaryInfC{$G'$}
\RightLabel{$\mathsf{(EW)}$}
\UnaryInfC{$G'\cup\{\Gamma',(\bigwedge_{i\in\omega}\phi_i)^*\rhd\Delta'\}$}
\end{prooftree}
by 
\begin{prooftree}
\AxiomC{$\vdots$}
\UnaryInfC{$\overline{G}'\cup G$}
\RightLabel{$\mathsf{(EW)}$}
\UnaryInfC{$\overline{G}'\cup G\cup\{\Gamma',\Gamma\rhd\Delta'\}$}
\end{prooftree}
\end{enumerate}
As $\phi$ is again not the cut-formula of any cut in $d_1$, one can verify the correctness of the resulting proof $d$ by transfinite induction on $\vert T\vert\oplus w(T)\oplus W(T)$ of subdervations $T$ of $d_1$.

Using the same reasoning as before, we also again derive $\rho(d)\leq\vert\phi\vert$ as well as $\vert T'\vert\leq 2\otimes (\vert d_0^{-1}\vert\oplus\vert T\vert)$ where $T'$ is the replacement of any rooted subderivation $T$ of $d_1$ \emph{after} the replacement and correction procedure, as before. The latter implies $\vert d\vert\leq 2\otimes (\vert d_0\vert\oplus\vert d_1\vert)$ as before.\\

The cases of $\phi=\phi_0\lor\phi_1$ and $\phi=\exists x\phi_0$ are similar to that of $\phi=\bigvee_{i\in\omega}\phi_i$ (see \cite{BC2002} for the latter). The quantifier case in particular uses Lemma \ref{lem:sub}.

Similarly, the cases of $\phi=\phi_0\land\phi_1$, $\phi=\forall x\phi_0$ and $\phi=\phi_0\rightarrow\phi_1$ are analogous to that of $\phi=\bigwedge_{i\in\omega}\phi_i$ (see again \cite{BC2002} for the last). As with $\exists$, the $\forall$-case uses Lemma \ref{lem:sub}.
\end{proof}
\begin{lemma}\label{lem:cutelim}
Let $d\vdash^s H$ with $\rho(d)\leq v+\omega^z$. Then there is a derivation $d'\vdash^s H$ with $\vert d'\vert\leq\chi^z(\vert d\vert)$ and $\rho(d')\leq v$.
\end{lemma}
\begin{proof}
The theorem is proved by induction on the lexicographically ordered pair $(z,\vert d\vert)$. So assume the claim for any $(\hat{z},\vert\hat{d}\vert)$ with $\hat{z}<z$ or $\hat{z}=z$ and $\vert\hat{d}\vert<\vert d\vert$. We divide the proof on whether the last inference rule of $d$ was a cut, a weakening or neither.\\

Suppose the last rule was not a cut and not a weakening. Let $k\leq\omega$ be the arity of the last rule and let $d_i$, $i<k$, be the direct predecessors with $d_i\vdash^sH_i$. Naturally, $\vert d_i\vert <\vert d\vert$ for all $i<k$ as the last rule was not a weakening and also
\[
v+\omega^z\geq\rho(d)\geq\rho(d_i)
\]
for all $i<k$ by definition. Using the induction hypothesis on $d_i$, we get derivations $d_i'\vdash^sH_i$ with $\vert d_i'\vert\leq\chi^z(\vert d_i\vert)$ and $\rho(d_i')\leq v$. Using the same last rule as in $d$, we combine the $d_i'$ to a proof $d'\vdash H$. First, we have $\rho(d')\leq v$ as the last rule was not a cut. Further, we get
\[
\vert d_i'\vert\leq\chi^z(\vert d_i\vert)<\chi^z(\vert d\vert)
\]
using that $\vert d_i\vert<\vert d\vert$ and that $\chi^z$ is increasing. As $\vert d'\vert$ is the least ordinal $\alpha$ with $\vert d'_i\vert <\alpha$ for all $i$, as the last rule was not a weakening, we get $\vert d'\vert\leq\chi^z(\vert d\vert)$.\\

Suppose the last rule was a cut. Then we get two preceding derivations
\[
d_0\vdash^s G\cup\{\Gamma\rhd\phi\}\text{ and }d_1\vdash^s G\cup\{\Gamma,\phi\rhd\Delta\}
\]
and by definition, we have
\[
v+\omega^z\geq\rho(d)=\max\{\vert\phi\vert+1,\rho(d_0),\rho(d_1)\}.
\]
and, as the last rule was not a weakening, we get $\vert d_i\vert <\vert d\vert$. We can apply the induction hypothesis to $d_0,d_1$ to get derivations $d'_0\vdash^s G\cup\{\Gamma\rhd\phi\}$ and $d_1'\vdash^s G\cup\{\Gamma,\phi\rhd\Delta\}$ with $\vert d_i'\vert\leq\chi^z(\vert d_i\vert)$ and $\rho(d_i')\leq v$.

If $z=0$, then $\rho(d)\leq v+\omega^z=v+1$ and therefore $\vert\phi\vert\leq v$.  Now, either
\begin{enumerate}
\item[(i)] $\max\{\rho(d_0'),\rho(d_1')\}\leq\vert\phi\vert$, or 
\item[(ii)] $\max\{\rho(d_0'),\rho(d_1')\}>\vert\phi\vert$.
\end{enumerate}
For (i), we may apply Lemma \ref{lem:reduction} to get a derivation
\[
d'\vdash G\cup\{\Gamma\rhd\Delta\}
\]
with $\rho(d')\leq \vert\phi\vert\leq v$ and
\begin{align*}
\vert d'\vert&\leq 2\otimes(\chi^0(\vert d_0\vert)\oplus\chi^0(\vert d_1\vert))\\
&\leq 2\otimes(4^{\otimes\vert d_0\vert}\oplus 4^{\otimes\vert d_1\vert})\\
&\leq 2\otimes (2\otimes 4^{\otimes\max\{\vert d_0\vert,\vert d_1\vert\}})\\
&\leq 4^{\otimes(\max\{\vert d_0\vert,\vert d_1\vert\}+1)}\\
&=4^{\otimes\vert d\vert}=\chi^0(\vert d\vert)
\end{align*}
where we in particular used that $\oplus$ and $\otimes$ are increasing in both arguments as well as commutative, associative and distributive, that $\chi^0$ is increasing, that $\alpha\oplus 1=\alpha+1$,  and that $\max\{\vert d_0\vert,\vert d_1\vert\}+1=\vert d\vert$ by definition.

For (ii), we combine $d_0'$ and $d_1'$ using cut on $\phi$ to a derivation $d'\vdash G\cup\{\Gamma\rhd\Delta\}$. Now, we get
\[
\rho(d)=\max\{\vert\phi\vert+1,\rho(d_0'),\rho(d_1')\}=\max\{\rho(d_0'),\rho(d_1')\}\leq v
\]
as $\vert\phi\vert+1\leq\max\{\rho(d_0'),\rho(d_1')\}$ by the assumption (ii) and as $\rho(d_i')\leq v$ from the induction hypothesis. Further, we have
\[
\vert d'\vert=\max\{\vert d_0'\vert,\vert d_1'\vert\}+1\leq \max\{4^{\otimes\vert d_0\vert},4^{\otimes\vert d_1\vert}\}+1\leq 4^{\otimes\vert d\vert}=\chi^0(\vert d\vert)
\]
which completes the case for (ii).\\

If $z\neq 0$, then there are $y<z$ and $k\in\mathbb{N}$ such that $\vert\phi\vert\leq v+\omega^y\cdot k$. We can combine $d_0'$ and $d_1'$ using cut to a derivation $\hat{d}\vdash^s G\cup\{\Gamma\rhd\Delta\}$.

This immediately gives $\vert\hat{d}\vert\leq\max\{\chi^z(\vert d_0\vert),\chi^z(\vert d_1\vert)\}+1$ and $\rho(\hat{d})\leq v+\omega^y\cdot k$. As $y<z$, we can apply the induction hypothesis $k$-times and get a derivation $d'\vdash^s G\cup\{\Gamma\rhd\Delta\}$ with $\rho(d')\leq v$ as well as
\[
\vert d'\vert\leq (\chi^y)^{(k)}(\max\{\chi^z(\vert d_0\vert),\chi^z(\vert d_1\vert)\}+1)
\]
We have $\vert d_i\vert <\vert d\vert$ and thus $\chi^z(\vert d_i\vert)<\chi^z(\vert d\vert)$ which implies $\chi^z(\vert d_i\vert)+1\leq\chi^z(\vert d\vert)$. This gives $\max\{\chi^z(\vert d_0\vert),\chi^z(\vert d_1\vert)\}+1\leq\chi^z(\vert d\vert)$. 

Now, assuming $\alpha\leq\chi^z(\beta)$ and $y<z$, then
\[
\chi^y(\alpha)\leq\chi^y(\chi^z(\beta))=\chi^z(\beta)
\]
where the first inequality follows from $\chi^y$ being increasing and the latter equality follows from the definition of $\chi^z(\beta)$ being the $\beta$-th simultaneous solution for $\gamma=\chi^x(\gamma)$ for all $x<z$. This implies, in combination with the above, that
\[
\vert d'\vert\leq (\chi^y)^{(k)}(\max\{\chi^z(\vert d_0\vert),\chi^z(\vert d_1\vert)\}+1)\leq\chi^z(\vert d\vert).
\]
\end{proof}
From Lemma \ref{lem:cutelim}, cut elimination immediately follows:
\begin{theorem}[Cut elimination]
For any derivation $d\vdash H$, there exists a cut-free derivation $d'\vdash H$.
\end{theorem}
Further, we of course get a bound on $\vert d'\vert$ which, in the case of finitary proofs with finitary formulas, matches that of Baaz and Ciabattoni \cite{BC2002}.
\section{The Range of the Results and Extensions}
We want to use this section to give an overview of some other topics extending the previous ones and some initial observations, at varying depth, regarding those.\\

\subsection{Extensions of the Completeness Results}
At first, it should be noted that the completeness theorems for $\mathsf{G}_{\omega_1}$ and $\mathsf{G}_{\omega_1,\omega}$ do not generalize to uncountable sets $\Gamma$. Consider the following two notions of compactness from \cite{Agu2016}:
\begin{definition}[essentially Aguilera \cite{Agu2016}]
Let $\Vdash\subseteq\mathcal{P}(\mathcal{L}_{\kappa,\lambda})\times\mathcal{L}_{\kappa,\lambda}$ or $\Vdash\subseteq\mathcal{P}(\mathcal{L}_{\kappa})\times\mathcal{L}_{\kappa}$ be a relation. Then $\Vdash$ is called
\begin{enumerate}
\item \emph{weakly compact} if for every $\Gamma,\phi$ with at most $\kappa$ many different atomics, there is a $\Delta$ with size $<\kappa$ such that $\Gamma\Vdash\phi$ implies $\Delta\Vdash\phi$,
\item \emph{compact} if for every $\Gamma,\phi$, there is a $\Delta$ with size $<\kappa$ such that $\Gamma\Vdash\phi$ implies $\Delta\Vdash\phi$.
\end{enumerate}
\end{definition}
Indeed, the relations $\vdash_{\mathcal{G}^{(D)}_{\omega_1}}$ and $\vdash_{\mathcal{G}^{(D)}_{\omega_1,\omega}}$ are compact as any proof only involves countably many formulas.

However, as Aguilera shows in \cite{Agu2016}, the consequence relation $\models_{\mathsf{G}_{\omega_1,\omega}}$ (i.e. $\models_{\mathbf{[0,1]_\mathbb{R}}}$) is not even weakly compact. This follows from the following general result:
\begin{proposition}[Aguilera \cite{Agu2016}]\label{pro:compact}
If $\models_{\mathsf{G}_{\kappa,\omega}}$ is weakly compact, then $\kappa$ is weakly compact, i.e. it is strongly inaccessible and any tree of size $\kappa$, such that every level has $<\kappa$ nodes, has a branch of length $\kappa$.
\end{proposition}
As $\omega_1$ is not strongly inaccessible, it is not weakly compact and therefore $\models_{\mathsf{G}_{\omega_1,\omega}}$ is not either.\\

It turns out that Aguilera's proof is, in itself, ``propositional" and can be straightforwardly adapted to $\mathsf{G}_{\kappa}$. We therefore have the following result which implies the same limitation for $\mathcal{G}^D_{\omega_1}$:
\begin{proposition}
If $\models_{\mathsf{G}_{\kappa}}$ is weakly compact, then $\kappa$ is weakly compact.
\end{proposition}
We omit the proof as it is, as hinted above, literally that of Aguilera \cite{Agu2016}, rephrased in the propositional language.

Before moving on to other topics, we want to mention a peculiar application of the above result of Aguilera and its propositional version. Although we don't dive into the rich (and difficult) topic of interpolation for G\"odel logics, we can use the impossibility of a proof calculus for uncountable premises to show the following negative result.
\begin{proposition}
There is no countable $\Delta\subseteq\mathcal{L}_{\omega_1,\omega}$ such that for any $\phi,\psi\in\mathcal{L}_{\omega_1,\omega}$, $\models_{\mathsf{G}_{\omega_1,\omega}}\phi\rightarrow\psi$ implies that there is a $\delta\in\Delta$ with $\models_{\mathsf{G}_{\omega_1,\omega}}\phi\rightarrow\delta$ and $\models_{\mathsf{G}_{\omega_1,\omega}}\delta\rightarrow\psi$. Similarly for $\mathcal{L}_{\omega_1}$ and $\mathsf{G}_{\omega_1}$.
\end{proposition}
\begin{proof}
Suppose such a $\Delta$ would exist. Constructing the Lindenbaum-Tarski algebra over $\mathcal{L}_{\omega_1,\omega}$, the set $\{[\delta]\mid\delta\in\Delta\}$ would form a countable dense subset of $\mathbf{LT}^\Gamma$, even for \emph{uncountable} $\Gamma$. As $\mathbf{LT}^\Gamma$ is therefore separable, it embeds into $[0,1]_\mathbb{R}$ with an embedding preserving infima and suprema. If we assume $\Gamma\not\vdash_{\mathcal{G}^D_{\omega_1,\omega}}\phi$ with said uncountable $\Gamma$, then this embedding would provide a countermodel by Lemma \ref{lem:fohomprop}, verifying $\Gamma\not\models_{\mathsf{G}_{\omega_1,\omega}}\phi$. Thus, we would have completeness of $\mathcal{G}^D_{\omega_1,\omega}$ for $\mathsf{G}_{\omega_1,\omega}$ w.r.t. uncountable $\Gamma$ which is a contradiction to Proposition \ref{pro:compact}.

The argument works similarly also for the propositional $\mathsf{G}_{\omega_1}$.
\end{proof}
Although it is probably expected that there is no countable set of interpolants in this infinitary case (as we already work over an uncountable language), it is maybe still instructive to note how cardinality considerations can have an impact on these type of questions.\\

In a similar vein, a generalization of the results to $\mathcal{L}_{\omega_1,\omega_1}$ is problematic. As is well-known by a theorem of Scott (and Karp, see \cite{Kar1964}\footnote{As Karp remarks, the proof of hers is based on an outline circulated by Scott in 1960 which was not published.}), the set of classical validities over $\mathcal{L}_{\omega_1,\omega_1}$ is not definable in $H(\omega_1)$, the collection of hereditarily countable sets, and thus in particular not $\Sigma_1$ on $H(\omega_1)$, a property which would, however, be implied by the existence of a complete classical proof calculus with proofs of countable lengths. These results should generalize to the G\"odel case:
\begin{question}
Is the set of theorems of $\mathsf{G}_{\omega_1,\omega_1}$ non-definable over $H(\omega_1)$?
\end{question}
We also want to note that there is a different definition of consequence common in the context of G\"odel logics, which we may define by
\[
\Gamma\models^\leq_{\mathsf{G}_{\kappa}}\phi\text{ if }\inf v[\Gamma]\leq v(\phi)\text{ for any evaluation }v
\]
for the propositional infinitary case and similarly in the first-order case by
\[
\Gamma\models^\leq_{\mathsf{G}_{\kappa,\lambda}}\phi\text{ if }\inf \overline{\mathfrak{I}}[\Gamma]\leq \overline{\mathfrak{I}}(\phi)\text{ for any interpretation }\mathfrak{I}.
\]
In a finitary context, it can be easily seen that $\models^\leq$ is equivalent to $\models$ and we can show a similar statement here if we restrict to countable sets. The following results are thereby natural generalizations of the finitary cases as given in \cite{BPZ2007} and the proofs for both results are essentially the same (and thus omitted).
\begin{lemma}\label{lem:cutoffeval}
Let $\mathbf{A}$ be a complete linear Heyting algebra, $x\in A$, and let $v:\mathcal{L}_{\kappa}\to\mathbf{A}$ be any evaluation. Then, for
\[
v_x(p):=\begin{cases}v(p)&\text{if }v(p)<^\mathbf{A}x,\\1^\mathbf{A}&\text{otherwise},\end{cases}
\]
with $p\in Var_{\kappa}$, the unique extension $\overline{v_x}:\mathcal{L}_{\kappa}\to\mathbf{A}$ to an evaluation satisfies:
\[
x\not\in v[\mathrm{sub(\phi)}]\text{ implies }\overline{v_x}(\phi)=\begin{cases}v(\phi)&\text{if }v(\phi)<^\mathbf{A}x,\\1^\mathbf{A}&\text{otherwise},\end{cases}
\]
for any $\phi\in\mathcal{L}_{\kappa}$.
\end{lemma}
The following first-order version is rather intricate to formulate if one wants to guarantee a similar level of generality as in the finitary case. We give it here in its full strength but will, in the following, mostly use the special case with $\kappa=\omega_1$ and $\lambda=\omega$ where the conditions simplify considerably (in particular the definition of $\mathrm{Val}_\mathfrak{I}(\phi)$).
\begin{lemma}\label{lem:cutoffevalfo}
Let $\mathbf{A}$ be a complete linear Heyting algebra, $x\in A$, and let $\mathfrak{M}$ be an $\mathbf{A}$-valued model. Define $\mathfrak{M}_x$ from $\mathfrak{M}$ by replacing $P^\mathfrak{M}$ with
\[
P^{\mathfrak{M}_x}(m_1,\dots,m_n):=\begin{cases}P^\mathfrak{M}(m_1,\dots,m_n)&\text{if }P^\mathfrak{M}(m_1,\dots,m_n)<^\mathbf{A}x,\\1^\mathbf{A}&\text{otherwise},\end{cases}
\]
where $P$ is an $n$-ary predicate. We write $\mathfrak{J}_x=(\mathfrak{M}_x,w)$ and
\begin{align*}
&\mathrm{Val}_{\mathfrak{J}}(\phi):=\Big\{\overline{\mathfrak{J}'}(\psi)\mid \psi\in\mathrm{sub}(\phi)\text{ and }\mathfrak{J}'=\left(\dots\left(\mathfrak{J}\tfrac{f_1}{X_1}\right)\tfrac{f_2}{X_2}\dots\right)\tfrac{f_n}{X_n}\\&\qquad\qquad\qquad\text{ where } X_i\subseteq Var_\kappa\cap\mathrm{sub}(\phi)\text{ with }\vert X_i\vert <\lambda\text{ and }f_i:X_i\to M\Big\}.
\end{align*}
given an interpretation $\mathfrak{J}=(\mathfrak{M},w)$. 

Then, for a given $v:Var_\kappa\to M$ and $\mathfrak{I}=(\mathfrak{M},v)$, we have
\[
x\not\in\mathrm{Val}_{\mathfrak{I}'}(\phi)\text{ implies }\overline{\mathfrak{I}'_x}(\phi)=\begin{cases}\overline{\mathfrak{I}'}(\phi)&\text{if }\overline{\mathfrak{I}'}(\phi)<^\mathbf{A}x,\\1^\mathbf{A}&\text{otherwise},\end{cases}
\]
for any $\phi\in\mathcal{L}_{\kappa,\lambda}$ and any interpretation
\[
\mathfrak{I}'=\left(\dots\left(\mathfrak{I}\tfrac{f_1}{X_1}\right)\tfrac{f_2}{X_2}\dots\right)\tfrac{f_n}{X_n}
\]
where $X_i\subseteq Var_\kappa$ with $\vert X_i\vert <\lambda$ and $f_i:X_i\to M$.
\end{lemma}
As a direct consequence, we obtain the following result:
\begin{proposition}
For any countable $\Gamma\cup\{\phi\}\subseteq\mathcal{L}_{\omega_1}$, $\Gamma\models^\leq_{\mathsf{G}_{\omega_1}}\phi$ iff $\Gamma\models_{\mathsf{G}_{\omega_1}}\phi$. Similarly, for any countable $\Gamma\cup\{\phi\}\subseteq\mathcal{L}_{\omega_1,\omega}$, we have $\Gamma\models^\leq_{\mathsf{G}_{\omega_1,\omega}}\phi$ iff $\Gamma\models_{\mathsf{G}_{\omega_1,\omega}}\phi$.
\end{proposition}
\subsection{Other Sets of Variables}
As common in propositional and first-order G\"odel logics, one could consider closed sets $V$ with $\{0,1\}\subseteq V\subsetneq [0,1]$ instead of $[0,1]$ as truth-value sets, thereby forming the propositional variants $\mathsf{G}^V_\kappa$ and the first-order variants $\mathsf{G}^V_{\kappa,\lambda}$ by extending the semantic definitions from before. The most common instances of $V$ are among
\begin{align*}
&V_\mathbb{R}:=[0,1],\\
&V_0:=\{0\}\cup[1/2,1],\\
&V_\downarrow:=\{1/k\mid k\geq 1\}\cup\{0\},\\
&V_\uparrow:=\{1-1/k\mid k\geq 1\}\cup\{1\},\\
&V_n:=\{1-1/k\mid 1\leq k\leq n-1\}\cup\{1\}\text{ with }n\geq 2,
\end{align*}
following the selection from \cite{BPZ2007}. In that notation, we have $\mathsf{G}_\kappa=\mathsf{G}^{V_\mathbb{R}}_\kappa$ and $\mathsf{G}_{\kappa,\lambda}=\mathsf{G}^{V_\mathbb{R}}_{\kappa,\lambda}$. A few easy observations from \cite{BPZ2007} directly carry over to the infinitary case.
\begin{proposition}\label{pro:relationships}
Let $\mathsf{G}^V$ be either $\mathsf{G}^V_\kappa$ or $\mathsf{G}^V_{\kappa,\lambda}$ for arbitrary $\kappa,\lambda$ (for $\lambda\leq\kappa$). The following relations hold:
\begin{enumerate}
\item $\mathsf{G}^{V_\mathbb{R}}=\bigcap_{V}\mathsf{G}^V$
\item $\mathsf{G}^{V_n}\supsetneq\mathsf{G}^{V_{n+1}}$
\item $\mathsf{G}^{V_n}\supsetneq\mathsf{G}^{V_{\uparrow}}\supsetneq\mathsf{G}^{V_\mathbb{R}}$
\item $\mathsf{G}^{V_n}\supsetneq\mathsf{G}^{V_{\downarrow}}\supsetneq\mathsf{G}^{V_\mathbb{R}}$
\item $\mathsf{G}^{V_{0}}\supsetneq\mathsf{G}^{V_\mathbb{R}}$.
\end{enumerate}
\end{proposition}
We omit the proof as it is essentially a replica of the analogous result in the finitary first-order case (see \cite{BPZ2007}). Still, we want to emphasize the following differences to the finitary case: already in the finitary propositional case, $\mathsf{G}^{V_\uparrow}_\omega$, $\mathsf{G}^{V_\downarrow}_\omega$ and $\mathsf{G}^{V_0}_\omega$ differ in entailment (see \cite{BZ1998}) while they coincide, as observed first by Dummett \cite{Dum1959}, in tautologies unlike the finitary first-order versions. So, it is natural that in this infinitary case the propositional (and first-order) versions differ in entailment as well.

But, in the infinitary case, the propositional variants already differ with respect to tautologies and, moreover, the witnessing (non-)tautologies are natural analogues of the finitary first-order examples: consider $C^\uparrow:=\bigvee_{i\in\omega}(p_i\rightarrow\bigwedge_{j\in\omega}p_j)$, $C^\downarrow:=\bigvee_{i\in\omega}p_i(\bigvee_{j\in\omega}p_j\rightarrow p_i)$ and $\mathsf{ISO}_0:=\bigwedge_{i\in\omega}\neg\neg p_i\rightarrow\neg\neg\bigwedge_{i\in\omega}p_i$. 

Then $C^\uparrow$ is valid in $\mathsf{G}^{V_\uparrow}$ but not in $\mathsf{G}^{V_{\downarrow}}$.  $C^\downarrow$ is valid in $\mathsf{G}^{V_\uparrow}$ and $\mathsf{G}^{V_{\downarrow}}$. Both are not valid in $\mathsf{G}^{V_0}$ and $\mathsf{G}^{V_\mathbb{R}}$. $\mathsf{ISO}_0$ is valid in $\mathsf{G}^{V_0}$ but not in $\mathsf{G}^{V_\mathbb{R}}$.\\

Further, we can give the following analogy of the relationship between $\mathsf{G}^{V_\uparrow}$ and the finite-valued $\mathsf{G}^{V_n}$.
\begin{proposition}\label{pro:vupequivcapfin}
We have
\[
\mathsf{G}^{V_\uparrow}_{\kappa}=\bigcap_{n\geq 2}\mathsf{G}^{V_n}_{\kappa}\text{ and }\mathsf{G}^{V_\uparrow}_{\kappa,\lambda}=\bigcap_{n\geq 2}\mathsf{G}^{V_n}_{\kappa,\lambda}.
\]
\end{proposition}
\begin{proof}
Again, let $\mathsf{G}^V$ be either $\mathsf{G}^V_{\kappa}$ or $\mathsf{G}^V_{\kappa,\lambda}$. Item (3) of Proposition \ref{pro:relationships} gives
\[
\mathsf{G}^{V_\uparrow}\subseteq\bigcap_{n\geq 2}\mathsf{G}^{V_n}.
\]
For the converse, suppose that $\Gamma\not\models_{\mathsf{G}^{V_\uparrow}}\phi$, i.e. there is an evaluation $v$ such that $v[\Gamma]\subseteq\{1\}\text{ but }v(\phi)<1$ in the propositional case or an interpretation $\mathfrak{I}$ such that $\overline{\mathfrak{I}}[\Gamma]\subseteq\{1\}\text{ but }\overline{\mathfrak{I}}(\phi)<1$ in the first-order case. As $v$/$\overline{\mathfrak{I}}$ evaluate into $V_\uparrow$, there is a $k$ such that $v(\phi)=1-\tfrac{1}{k}$ or $\overline{\mathfrak{I}}(\phi)=1-\tfrac{1}{k}$. Let $x\in [0,1]$ be such that $1-\tfrac{1}{k}<x<1-\tfrac{1}{k+1}$ and
\[
x\not\in v[\mathrm{sub}(\Gamma\cup\{\phi\})]
\]
in the propositional case or such that
\[
x\not\in \mathrm{Val}_{\mathfrak{I}}(\Gamma\cup\{\phi\})
\]
in the first-order case. We form $v_x$ or $\mathfrak{I}_x$ by Lemma \ref{lem:cutoffeval} or Lemma \ref{lem:cutoffevalfo}, respectively. The above choice of $x$ is such that $\overline{v_x}[\Gamma]\subseteq\{1\}\text{ but }\overline{v_x}(\phi)<1$ in the propositional case and $\overline{\mathfrak{I}_x}[\Gamma]\subseteq\{1\}\text{ but }\overline{\mathfrak{I}_x}(\phi)<1$ in the first-order case by the previous lemmas. But, by the choice of $x$, we have that $\overline{v_x}$ or $\overline{\mathfrak{I}_x}$ evaluate into $V_{k+1}$ which gives $\Gamma\not\models_{\mathsf{G}^{V_{k+1}}}\phi$.
\end{proof}
By the results of \cite{BPZ2007,Dum1959}, the status quo on complete proof calculi in the finitary setting in very clear cut: $\mathsf{G}^V_{\omega,\omega}$ is axiomatizable iff $V$ is finite or uncountable with either $0$ contained in the perfect kernel of $V$ or isolated. In particular, already the tautologies of $\mathsf{G}^{V_\uparrow}_{\omega,\omega}$ and $\mathsf{G}^{V_\downarrow}_{\omega,\omega}$ are not recursively enumerable. On the propositional side, while the tautologies of all $\mathsf{G}^V_\omega$ are axiomatizable (again, see\cite{Dum1959}), the only axiomatizable entailment relations are $\mathsf{G}^{V_n}_\omega$ and $\mathsf{G}^{V_\mathbb{R}}_\omega$ (see \cite{BZ1998}).\\

Now, the situation is different in the infinitary cases. In the following, we will obtain analogous axiomatizations for the instances of $V$ which were axiomatizable already in the finite but we further obtain infinitary axiomatizations of $\mathsf{G}^{V_\uparrow}_{\omega_1,\omega}$ and $\mathsf{G}^{V_\uparrow}_{\omega_1}$.

We don't now the state of $\mathsf{G}^{V_\downarrow}_{\omega_1,\omega}$, $\mathsf{G}^{V_\downarrow}_{\omega_1}$ or any other $V$ in particular but in the finitary, as shown by H\'ajek \cite{Haj2005}, the tautologies of $\mathsf{G}^{V_\downarrow}_{\omega,\omega}$ are not arithmetical. So we leave with the following question regarding the other truth-value sets:
\begin{question}
Have any $\mathsf{G}^V_{\omega_1}$ or $\mathsf{G}^V_{\omega_1,\omega}$ (countable) infinitary axiomatizations for any $V$ not considered here?
\end{question}
To approach these axiomatizability questions, we follow the general route of \cite{BPZ2007} which relies on tools from the theory of Polish spaces, like the Cantor-Bendixson theorem, which we briefly want to recall. In the following, we write $\mathbf{A}_V$ for the Heyting algebra associated with a G\"odel set $V$. Note that every $V$, as a closed subset of $\mathbb{R}$, is a Polish space. A subset $P$ of $\mathbb{R}$ is \emph{perfect} if it is closed and every point is a limit point in the topology induced by $\mathbb{R}$.
\begin{theorem}[Cantor-Bendixson]
Any Polish space $X$ can be partitioned as $P\cup C$ with $P$ perfect and $C$ countable and open.
\end{theorem}
The following result is then the central connection between the Cantor-Bendixson theorem and evaluations over G\"odel sets.
\begin{lemma}[Preining \cite{Pre2003}]\label{lem:perfembed}
Let $M\subseteq[0,1]$ be countable and $P\subseteq [0,1]$ be perfect. Then there is a strictly monotone $h:M\to P$ which preserves any infima and suprema existing in $M$ and if $\inf M\in M$, then $h(\inf M)=\inf P$.
\end{lemma}
\subsubsection{$V$ is finite}
We consider the axiom scheme
\[
\mathsf{FIN}(n):=(\phi_0\rightarrow\phi_1)\lor(\phi_1\rightarrow\phi_2)\lor\dots\lor(\phi_{n-1}\rightarrow\phi_n)
\]
as in the finitary axiomatizations.
\begin{theorem}\label{thm:propfincomp}
For any countable $\Gamma\cup\{\phi\}\subseteq\mathcal{L}_{\omega_1}$, we have
\[
\Gamma\vdash_{\mathcal{G}^D_{\omega_1}+\mathsf{FIN}(n)}\phi\text{ iff }\Gamma\models_{\mathsf{G}^{V_n}_{\omega_1}}\phi.
\]
Similarly, for any countable $\Gamma\cup\{\phi\}\subseteq\mathcal{L}_{\omega_1,\omega}$ where all formulas of $\Gamma$ are closed, we have
\[
\Gamma\vdash_{\mathcal{G}^D_{\omega_1,\omega}+\mathsf{FIN}(n)}\phi\text{ iff }\Gamma\models_{\mathsf{G}^{V_n}_{\omega_1,\omega}}\phi.
\]
In fact, $\mathsf{FIN}(n)$ can, in both cases, be replaced by $\mathsf{FIN}^a(n)$: all atomic instances of $\mathsf{FIN}(n)$.
\end{theorem}
\begin{proof}
Soundness is routine. For the converse, define $\mathcal{L}_A:=\mathrm{frag}(\Gamma\cup\{\phi\})$ and write $\mathcal{L}^a_A$ for the atomics of $\mathcal{L}_A$. We consider
\[
\Pi:=\{(\phi_0\rightarrow\phi_1)\lor(\phi_1\rightarrow\phi_2)\lor\dots\lor(\phi_{n-1}\rightarrow\phi_n)\mid\phi_i\in\mathcal{L}^a_A\}.
\]
$\Pi$ is countable as $\Gamma$, and therefore $\mathcal{L}_A$, is countable. Now suppose $\Gamma\not\vdash_{(\mathcal{G}^D_{\omega_1}+\mathsf{FIN}(n))(\mathcal{L}_A)}\phi$. Then clearly $\Gamma\cup\Pi\not\vdash_{\mathcal{G}^D_{\omega_1}(\mathcal{L}_A)}\phi$ and by strong completeness of $\mathcal{G}^D_{\omega_1}(\mathcal{L}_A)$, we have $\Gamma\cup\Pi\not\models_{[0,1]_\mathbb{R}(\mathcal{L}_A)}\phi$, i.e. there is an evaluation $v:\mathcal{L}_{A}\to [0,1]$ with
\[
v[\Gamma\cup\Pi]\subseteq\{1\}\text{ but }v(\phi)<1.
\]
Now, the set $v[\mathcal{L}^a_A]$ contains at most $n$ elements. If not, then there are formulas $\phi_0,\dots,\phi_n\in\mathcal{L}^a_A$ with $v(\phi_i)>v(\phi_{i+1})$. In that case, we have
\[
v((\phi_0\rightarrow\phi_1)\lor(\phi_1\rightarrow\phi_2)\lor\dots\lor(\phi_{n-1}\rightarrow\phi_n))<1
\]
which is a contradiction to $v[\Pi]\subseteq\{1\}$. Thus, we can write $v[\mathcal{L}^a_A]\subseteq\{0,v_1,\dots,v_{n-2},1\}$ with $v_i<v_{i+1}$. By induction on the structure of formulas, we also get $v[\mathcal{L}_A]$ $\subseteq$ $\{0,v_1,\dots,v_{n-2},1\}$. We define a function $h:v[\mathcal{L}_A]\to V_n$ by setting $h(0):=0$, $h(1):=1$ and $h(v_i):=1-\tfrac{1}{i+1}$. $v[\mathcal{L}_A]$ is, with its order by $<$, a Heyting algebra and therefore $h$ is an isomorphism of Heyting algebras and in particular preserves infima and suprema. Lemma \ref{lem:homprop} gives that $h\circ v$ is a $\mathcal{L}_A$-evaluation with
\[
(h\circ v)[\Gamma]\subseteq\{1\}\text{ but }(h\circ v)(\phi)<1.
\]
As $h\circ v$ evaluates into $V_n$, we get $\Gamma\not\models_{\mathbf{A}_{V_n}(\mathcal{L}_A)}\phi$. Therefore also $\Gamma\not\models_{\mathsf{G}^{V_n}_{\omega_1}}\phi$, as $\mathbf{A}_{V_n}$ is complete.\\

For the first-order instances, we instead consider
\[
\Pi:=\{\forall x_{i_1}\dots\forall x_{i_m}((\phi_0\rightarrow\phi_1)\lor\dots\lor(\phi_{n-1}\rightarrow\phi_n))\mid \phi_j\in\mathcal{L}^a_A, \mathrm{var}(\phi_j)\subseteq\{x_{i_1},\dots,x_{i_m}\}\}.
\]
similar to \cite{BPZ2007} where $\mathcal{L}_A^a$ now represents the atomic formulas of the first-order fragment $\mathcal{L}_A$.

For countable fragments $\mathcal{L}_A$, $\Pi$ is countable and as every atomic formulas has only a finite number of variables, we have that every formula in $\Pi$ is closed. Then apply Lemma \ref{lem:fohomprop} in place of Lemma \ref{lem:homprop} as in \cite{BPZ2007}.
\end{proof}
\subsubsection{$V$ is $V_\uparrow$}
The strength of infinitary logics is of course that we have infinitary disjunctions available which we can use to combine the various finitary axioms $\mathsf{FIN}(n)$. More precisely, we define the scheme $\mathsf{FIN}$ by
\[
\bigvee_{n\geq 2}\bigwedge_{k\in\omega}\bigvee_{i=0}^{n-1}\left(\phi^{n,k}_{i}\rightarrow\phi^{n,k}_{i+1}\right).
\]
In the first-order case, we will additionally consider a seemingly weakened version $\mathsf{FIN}^a$ given by
\[
\bigvee_{n\geq 2}\bigwedge_{k\in\omega}\bigvee_{i=0}^{n-1}\forall x^{n,k}_{j_1}\dots\forall x^{n,k}_{j_m}\left(\phi^{n,k}_{i}\rightarrow\phi^{n,k}_{i+1}\right)
\]
where all $\phi^{n,k}_i$ are atomic with $\mathrm{var}(\phi^{n,k}_i)\subseteq\{x^{n,k}_{j_1},\dots,x^{n,k}_{j_m}\}$. Now, $\mathsf{FIN}$ can be used to obtain an axiomatization of $\bigcap_{n\geq 2}\mathsf{G}^{V_n}_{\omega_1}$ or $\bigcap_{n\geq 2}\mathsf{G}^{V_n}_{\omega_1,\omega}$ for countable sets. In combination with the previous Proposition \ref{pro:vupequivcapfin}, we then obtain an axiomatization of $V_\uparrow$.
\begin{theorem}
For any countable $\Gamma\cup\{\phi\}\subseteq\mathcal{L}_{\omega_1}$, we have
\[
\Gamma\vdash_{\mathcal{G}^D_{\omega_1}+\mathsf{FIN}}\phi\text{ iff }\Gamma\models_{\mathsf{G}^{V_\uparrow}_{\omega_1}}\phi.
\]
Similarly, for countable $\Gamma\cup\{\phi\}\subseteq\mathcal{L}_{\omega_1,\omega}$ where all formulas in $\Gamma$ are closed, we have
\[
\Gamma\vdash_{\mathcal{G}^D_{\omega_1,\omega}+\mathsf{FIN}}\phi\text{ iff }\Gamma\models_{\mathsf{G}^{V_\uparrow}_{\omega_1,\omega}}\phi.
\]
\end{theorem}
\begin{proof}
$\mathsf{FIN}$ is valid in $V_\uparrow$: suppose $v:\mathcal{L}_{\omega_1}\to V_\uparrow$ is such that 
\[
v\left(\bigvee_{n\geq 2}\bigwedge_{k\in\omega}\bigvee_{i=0}^{n-1}\left(\phi^{n,k}_{i}\rightarrow\phi^{n,k}_{i+1}\right)\right)<1.
\]
Thus, for some $\alpha$, we have 
\[
v\left(\bigwedge_{k\in\omega}\bigvee_{i=0}^{n-1}\left(\phi^{n,k}_{i}\rightarrow\phi^{n,k}_{i+1}\right)\right)\leq\alpha<1
\]
for any $n\geq 2$. Let $k$ be such that $v\left(\bigvee_{i=0}^{n-1}\left(\phi^{n,k}_{i}\rightarrow\phi^{n,k}_{i+1}\right)\right)<1$. Such a $k$ exists as $\alpha<1$. For such a $k$, we get
\[
v(\phi_0^{n,k})>v(\phi_1^{n,k})>\dots>v(\phi^{n,k}_n)
\]
and as $v$ evaluates into $V_\uparrow$, we have
\[
v(\phi_0^{n,k})\geq 1-\tfrac{1}{n+1}\text{ and }v(\phi_1^{n,k})\geq 1-\tfrac{1}{n}
\]
and therefore $v(\phi_0^{n,k}\rightarrow\phi_1^{n,k})\geq 1-\tfrac{1}{n}$. This yields
\[
v\left(\bigvee_{i=0}^{n-1}\left(\phi^{n,k}_{i}\rightarrow\phi^{n,k}_{i+1}\right)\right)\geq 1-\tfrac{1}{n}
\]
for any such $k$ and we get
\begin{align*}
v\left(\bigwedge_{k\in\omega}\bigvee_{i=0}^{n-1}\psi^{n,k}\right)&=\min\left\{\inf_{v(\psi^{n,k})<1}v(\psi^{n,k}),\inf_{v(\psi^{n,k})=1}v(\psi^{n,k})\right\}\\
&=\inf_{v(\psi^{n,k})<1}v(\psi^{n,k})\\
&\geq 1-\tfrac{1}{n}
\end{align*}
for any $n\geq 2$ where we write $\psi^{n,k}:=\bigvee_{i=0}^{n-1}\left(\phi^{n,k}_{i}\rightarrow\phi^{n,k}_{i+1}\right)$. But this implies
\[
v\left(\bigvee_{n\geq 2}\bigwedge_{k\in\omega}\bigvee_{i=0}^{n-1}\left(\phi^{n,k}_{i}\rightarrow\phi^{n,k}_{i+1}\right)\right)=1
\]
in contradiction to our assumption.

For the converse, suppose $\Gamma\not\vdash_{\mathcal{G}^D_{\omega_1}+\mathsf{FIN}}\phi$. Then, we have
\[
\Gamma\not\vdash_{\mathcal{G}^D_{\omega_1}+\mathsf{FIN}(n)}\phi
\]
for some $n$ as if $\Gamma\vdash_{\mathcal{G}^D_{\omega_1}+\mathsf{FIN}(n)}\phi
$ for all $n$, then there are countably many instances
\[
\bigvee_{i=0}^{n-1}\left(\phi^{n,k}_{i}\rightarrow\phi^{n,k}_{i+1}\right) \quad(k\in\omega)
\]
of $\mathsf{FIN}(n)$, such that
\[
\Gamma\vdash_{\mathcal{G}^D_{\omega_1}}\bigwedge_{k\in\omega}\bigvee_{i=0}^{n-1}\left(\phi^{n,k}_{i}\rightarrow\phi^{n,k}_{i+1}\right)\rightarrow \phi
\]
for any $n$ and thus
\[
\Gamma\vdash_{\mathcal{G}^D_{\omega_1}}\bigvee_{n\geq 2}\bigwedge_{k\in\omega}\bigvee_{i=0}^{n-1}\left(\phi^{n,k}_{i}\rightarrow\phi^{n,k}_{i+1}\right)\rightarrow\phi
\]
by $(R\omega)_1$. The premise is an instance of $\mathsf{FIN}$ which implies that $\Gamma\vdash_{\mathcal{G}^D_{\omega_1}+\mathsf{FIN}}\phi$, a contradiction. Therefore, there is an $n$ with $\Gamma\not\vdash_{\mathcal{G}^D_{\omega_1}+\mathsf{FIN}(n)}\phi$ and thus, using that $\Gamma$ is countable, we get
\[
(\Gamma,\phi)\not\in\bigcap_{n\geq 2}\mathsf{G}^{V_n}_{\omega_1}
\]
by Theorem \ref{thm:propfincomp} which implies $\Gamma\not\models_{\mathsf{G}^{V_\uparrow}_{\omega_1}}\phi$ by Proposition \ref{pro:vupequivcapfin}.\\

The first order case is very similar. Soundness follows in the same way and the converse follows by the following slightly modified argument: supposing $\Gamma\not\vdash_{\mathcal{G}^D_{\omega_1,\omega}+\mathsf{FIN}}\phi$, we of course also have $\Gamma\not\vdash_{\mathcal{G}^D_{\omega_1,\omega}+\mathsf{FIN}^a}\phi$. As in the propositional case, we get $\Gamma\not\vdash_{\mathcal{G}^D_{\omega_1,\omega}+\mathsf{FIN}^a(n)}\phi$ for some $n$ where it is essential that all instances of $\mathsf{FIN}^a$ are closed to be able to use the deduction theorem. We get $\Gamma\not\models_{\mathsf{G}^{V_\uparrow}_{\omega_1,\omega}}\phi$ using Proposition \ref{pro:vupequivcapfin}.
\end{proof}
\subsubsection{$0$ contained in the perfect kernel}
We obtain the following infinitary version of the well-known finitary result from \cite{BPZ2007}.
\begin{lemma}\label{lem:perfkerequiv}
Let $V$ be a G\"odel set and $P$ its perfect Kernel and $W=V\cup[\inf P,1]$. For any fragment $\mathcal{L}_A$ (of $\mathcal{L}_{\omega_1}$ or $\mathcal{L}_{\omega_1,\omega}$) and any countable $\Gamma\cup\{\phi\}\subseteq\mathcal{L}_A$, we have
\[
\Gamma\models_{\mathbf{A}_V(\mathcal{L}_A)}\phi\text{ iff }\Gamma\models_{\mathbf{A}_{W}(\mathcal{L}_A)}\phi.
\]
\end{lemma}
\begin{proof}
We only present the propositional case which is essentially contained in \cite{BPZ2007}. The first order case is the same with Lemma \ref{lem:cutoffeval} replaced by Lemma \ref{lem:cutoffevalfo} and Lemma \ref{lem:homprop} replaced by Lemma \ref{lem:fohomprop}, respectively

Since $V\subseteq W$, we have that $\Gamma\models_{\mathbf{A}_W(\mathcal{L}_A)}\phi$ implies $\Gamma\models_{\mathbf{A}_{V}(\mathcal{L}_A)}\phi$. 

Let $v:\mathcal{L}_A\to\mathbf{A}_W$ be such that $v[\Gamma]\subseteq\{1\}$ but $v(\phi)<1$. As $\Gamma$ is countable, there is an $x\in [0,1]$ with $v(\phi)<x<1$ and where $x\not\in v[\mathrm{sub}(\Gamma\cup\{\phi\})]$. With $v_x$ as in Lemma \ref{lem:cutoffeval}, we have
\[
\overline{v_x}(\psi)=\begin{cases}v(\psi)&\text{if }v(\psi)<x,\\1&\text{otherwise},\end{cases}
\]
for any $\psi\in\Gamma\cup\{\phi\}$. Set $M:=\{v_x(\psi)\mid \psi\in\mathcal{L}_A\}\cup\{1\}$ and 
\[
M_0:=M\cap [0,\inf P)\text{ as well as }M_1:=(M\cap [\inf P,x])\cup\{\inf P\}.
\]
Lemma \ref{lem:perfembed} gives a strictly monotone $h:M_1\to P$ preserving all infima and suprema such that $h(\inf M_1)=\inf P$. We define
\[
g:y\mapsto \begin{cases}y&\text{if }y\in [0,\inf P],\\h(y)&\text{if }y\in [\inf P,x],\\1&\text{if }y=1,\end{cases}
\]
for $y\in M$. $g$ preserves infima, suprema and is strictly monotone with $g(0)=0$ and $g(1)=1$. $g$ therefore is a homomorphism of Heyting algebras preserving infima and suprema and by Lemma \ref{lem:homprop}, $g\circ v_x$ is an evaluation. But $g\circ v_x$ maps into $V$ and $(g\circ v_x)[\Gamma]\subseteq\{1\}$ but $g(v_x(\phi))=g(v(\phi))<1$.
\end{proof}
This immediately yields the following completeness result.
\begin{theorem}
Let $V$ be a G\"odel set where $0$ is contained in the perfect Kernel. For any fragment $\mathcal{L}_A$ (of $\mathcal{L}_{\omega_1}$ or $\mathcal{L}_{\omega_1,\omega}$) and countable $\Gamma\cup\{\phi\}\subseteq\mathcal{L}_A$ where $\Gamma$ is closed in the first-order case, we have
\[
\Gamma\vdash_{\mathcal{G}^D(\mathcal{L}_A)}\phi\text{ iff }\Gamma\models_{\mathbf{A}_V(\mathcal{L}_A)}\phi.
\]
Here, we write $\mathcal{G}^D(\mathcal{L}_A)$ for $\mathcal{G}^D_{\omega_1}(\mathcal{L}_A)$ or $\mathcal{G}^D_{\omega_1,\omega}(\mathcal{L}_A)$, respectively, depending the choice of language.
\end{theorem}
\subsubsection{$V$ uncountable and $0$ isolated}
We now turn to the case of an isolated $0$. For the first-order case, we have to additionally consider the quantifier version $\mathsf{QISO}_0$ of $\mathsf{ISO}_0$ (from which it was derived) as e.g. seen in \cite{BPZ2007}:
\[
\forall x\neg\neg\phi\rightarrow\neg\neg\forall x\phi.
\]
The following lemma is an easy adaption of the finitary first-order case from \cite{BPZ2007}, which only mentions the latter statement regarding the quantifiers.
\begin{lemma}\label{lem:isocons}
For any $\phi_i$, we have $\vdash_{\mathcal{G}^D_{\omega_1}+\mathsf{ISO_0}}\neg\bigwedge_{i\in\omega}\phi_i\rightarrow\bigvee_{i\in\omega}\neg\phi_i$. Similarly for $\mathcal{G}^D_{\omega_1,\omega}$. In the first-order case, we additionally have $\vdash_{\mathcal{G}^D_{\omega_1,\omega}+\mathsf{QISO}_0}\neg\forall x\phi\rightarrow\exists x\neg\phi$ for any $\phi$ and $x$.
\end{lemma}
\begin{theorem}\label{thm:comp0isocount}
Let $V$ be an uncountable G\"odel set with 0 isolated. Let $\mathcal{L}_A$ be a countable fragment of $\mathcal{L}_{\omega_1}$ or $\mathcal{L}_{\omega_1,\omega}$ where $\bigwedge_{i\in\omega}\phi_i\in\mathcal{L}_A$ iff $\bigvee_{i\in\omega}\phi_i\in\mathcal{L}_A$ and where $\bigvee_{i\in\omega}\phi_i\in\mathcal{L}_A$ implies $\bigvee_{i\in\omega}\neg\phi_i\in\mathcal{L}_A$. 

For any $\Gamma\cup\{\phi\}\subseteq\mathcal{L}_A$, we have
\[
\Gamma\vdash_{(\mathcal{G}^D_{\omega_1}+\mathsf{ISO_0})(\mathcal{L}_A)}\phi\text{ iff }\Gamma\models_{\mathbf{A}_V(\mathcal{L}_A)}\phi
\]
in the propositional case and
\[
\Gamma\vdash_{(\mathcal{G}^D_{\omega_1,\omega}+\mathsf{ISO_0}+\mathsf{QISO}_0)(\mathcal{L}_A)}\phi\text{ iff }\Gamma\models_{\mathbf{A}_V(\mathcal{L}_A)}\phi
\]
in the first-order case where $\Gamma$ is assumed to be closed.
\end{theorem}
\begin{proof}
Soundness is clear. For the other direction, note that we know by Lemma \ref{lem:perfkerequiv} that
\[
\Gamma\models_{\mathbf{A}_V(\mathcal{L}_A)}\phi\text{ iff }\Gamma\models_{\mathbf{A}_{V\cup[\inf P,1]}(\mathcal{L}_A)}\phi
\]
where $P$ is the perfect kernel of $V$. We thus assume that $[\inf P,1]\subseteq V$. We define
\[
\Pi:=\left\{\neg\bigwedge_{i\in\omega}\phi_i\rightarrow\bigvee_{i\in\omega}\neg\phi_i\mid \bigwedge_{i\in\omega}\phi_i\in\mathcal{L}_A\right\}.
\]
By the assumptions on $\mathcal{L}_A$, we have $\Pi\subseteq\mathcal{L}_A$ and in particular, $\Pi$ is countable. Suppose that $\Gamma\not\vdash_{(\mathcal{G}^D_{\omega_1}+\mathsf{ISO}_0)(\mathcal{L}_A)}\phi$.

We now either have $\Pi\cup\Gamma\models_{\mathbf{[0,1]_\mathbb{R}}(\mathcal{L}_A)}\phi$ or $\Pi\cup\Gamma\not\models_{\mathbf{[0,1]_\mathbb{R}}(\mathcal{L}_A)}\phi$. For the former, we get
\[
\Pi\cup\Gamma\vdash_{\mathcal{G}^D_{\omega_1}(\mathcal{L}_A)}\phi
\]
by completeness. As by Lemma \ref{lem:isocons}, $(\mathcal{G}^D_{\omega_1}+\mathsf{ISO}_0)(\mathcal{L}_A)$ proves every element of $\Pi$, we have
\[
\Gamma\vdash_{(\mathcal{G}^D_{\omega_1}+\mathsf{ISO}_0)(\mathcal{L}_A)}\phi
\]
which is a contradiction to our assumption. Therefore we have $\Pi\cup\Gamma\not\models_{\mathbf{[0,1]_\mathbb{R}}(\mathcal{L}_A)}\phi$, i.e. there is a $v:\mathcal{L}_{A}\to [0,1]$ with $v[\Pi]\cup v[\Gamma]\subseteq\{1\}$ and $v(\phi)<1$. We define
\[
h:x\mapsto\begin{cases}0&\text{if }x=0,\\\inf P+\frac{x}{1-\inf P}&\text{otherwise}.\end{cases}
\]
Note that $(h\circ v)[\mathcal{L}_A]\subseteq V$. Uniquely extend $v_h:x\mapsto h(v(x))$ for $x\in Var_{\omega_1}\cap \mathcal{L}_A$ to $\overline{v_h}$ on $\mathcal{L}_{A}$. Then we have
\[
\overline{v_h}(\psi)=h(v(\psi))
\]
for any $\psi\in\mathcal{L}_A$ which gives the claim. This can be proven by induction on $\psi$, see in particular \cite{BPZ2007} for the similar finitary case, where $v[\Pi]\subseteq\{1\}$ is used to handle the $\bigwedge$-case.

The proof in the first-order case is very similar. We then consider
\begin{align*}
&\Pi:=\left\{\forall\overline{x}\left(\neg\bigwedge_{i\in\omega}\phi_i\rightarrow\bigvee_{i\in\omega}\neg\phi_i\right)\mid \bigwedge_{i\in\omega}\phi_i\in\mathcal{L}_A,\overline{x}\in (Var_A)^n,n\in\mathbb{N}\right\}\\
&\qquad\qquad\qquad\cup\{\forall\overline{y}(\neg\forall x\phi\rightarrow\exists x\neg\phi)\mid \phi\in\mathcal{L}_A, x\in Var_A, \overline{y}\in (Var_A)^n,n\in\mathbb{N}\}
\end{align*}
and one proceeds as above and obtains a $\mathfrak{I}$ with $\overline{\mathfrak{I}}[\Pi]\cup\overline{\mathfrak{I}}[\Gamma]\subseteq\{1\}$ but $\overline{\mathfrak{I}}(\phi)<1$. As $\Pi$ is not closed, note in particular Remark \ref{rem:closednotneeded} in that closedness is not needed for both directions of the completeness results, only for the soundness direction. Note that $\Pi$ is countable.

With the same $h$ defined as before, one similarly defines $\mathfrak{I}'_h$ for $\mathfrak{I}'=\mathfrak{I}\tfrac{\overline{m}}{\overline{x}}$ by changing $P^\mathfrak{M}$ to $h\circ P^\mathfrak{M}$ for predicates $P$ in the underlying model $\mathfrak{M}$ and the key point is to now establish 
\[
\overline{\mathfrak{I}'_h}(\psi)=h(\overline{\mathfrak{I}'}(\psi))
\]
for any $\psi\in\mathcal{L}_A$ for any such $\mathfrak{I}'$ where, for the $\bigwedge$- and $\forall$-cases, it is important that $\overline{\mathfrak{I}}[\Pi]\subseteq\{1\}$ and that any element of $\Pi$ is universally quantified with arbitrary but finitely many quantifiers such that $\overline{\mathfrak{I}}[\Pi]\subseteq\{1\}$ implies 
\[
\overline{\mathfrak{I}'}\left(\neg\bigwedge_{i\in\omega}\phi_i\rightarrow\bigvee_{i\in\omega}\neg\phi_i\right)=1
\]
and 
\[
\overline{\mathfrak{I}'}(\neg\forall x\phi\rightarrow\exists x\neg\phi)=1
\]
even for any $\mathfrak{I}'$ as above.
\end{proof}
As before, we can now lift the above result to arbitrary $\mathcal{L}_A$ if we restrict to countably many assumptions.
\begin{corollary}
Let $V$ be an uncountable G\"odel set with 0 isolated and let $\mathcal{L}_A$ be an arbitrary fragment of $\mathcal{L}_{\omega_1}$ or $\mathcal{L}_{\omega_1,\omega}$ with the same additional closure conditions as in Theorem \ref{thm:comp0isocount} but let $\Gamma\cup\{\phi\}\subseteq\mathcal{L}_A$ be countable with $\Gamma$ closed in the first-order case.

For any $\Gamma\cup\{\phi\}\subseteq\mathcal{L}_A$, we have
\[
\Gamma\vdash_{(\mathcal{G}^D_{\omega_1}+\mathsf{ISO_0})(\mathcal{L}_A)}\phi\text{ iff }\Gamma\models_{\mathbf{A}_V(\mathcal{L}_A)}\phi
\]
in the propositional case and
\[
\Gamma\vdash_{(\mathcal{G}^D_{\omega_1,\omega}+\mathsf{ISO_0}+\mathsf{QISO}_0)(\mathcal{L}_A)}\phi\text{ iff }\Gamma\models_{\mathbf{A}_V(\mathcal{L}_A)}\phi
\]
in the first-order case.
\end{corollary}
The proof follows the same type of argument as in Corollary \ref{cor:fullcompgld} where we now consider the smallest fragment containing $\Gamma\cup\{\phi\}$ with the additional closure properties. The important point is here that for a countable $\Gamma$, this fragment is again countable.\\

Again, the above gives in particular
\[
\Gamma\vdash_{\mathcal{G}^D_{\omega_1}+\mathsf{ISO}_0}\phi\text{ iff }\Gamma\models_{\mathsf{G}^V_{\omega_1}}\phi
\]
for countable $\Gamma$ and
\[
\Gamma\vdash_{\mathcal{G}^D_{\omega_1,\omega}+\mathsf{ISO}_0+\mathsf{QISO}_0}\phi\text{ iff }\Gamma\models_{\mathsf{G}^V_{\omega_1,\omega}}\phi
\] 
for countable and closed $\Gamma$ where $0$ is isolated in $V$.
\subsection*{Acknowledgments}
I want to thank Matthias Baaz for helpful discussions of the topics of this paper.

\bibliographystyle{plain}
\bibliography{ref}

\end{document}